\documentclass[11pt,reqno]{amsart}
\usepackage[T1]{fontenc}
\usepackage{lmodern,microtype}
\usepackage[a4paper,margin=1in]{geometry}
\usepackage{mathtools,amssymb,amsthm,enumitem}
\usepackage[colorlinks=true,linkcolor=blue,citecolor=blue,urlcolor=blue]{hyperref}
\numberwithin{equation}{section}
\setlist[enumerate]{label=\textup{(\roman*)},leftmargin=*,itemsep=4pt,beginpenalty=10000}
\newcommand{\C}{\mathbb C}
\newcommand{\Z}{\mathbb Z}
\newcommand{\g}{\mathfrak g}
\newcommand{\s}{\mathfrak{sl}_2}
\newcommand{\T}{\mathcal T}
\DeclareMathOperator{\Der}{Der}
\DeclareMathOperator{\Hom}{Hom}
\DeclareMathOperator{\End}{End}
\DeclareMathOperator{\soc}{soc}
\DeclareMathOperator{\rad}{rad}
\DeclareMathOperator{\im}{im}
\DeclareMathOperator{\rank}{rank}
\DeclareMathOperator{\ad}{ad}
\newtheorem{theorem}{Theorem}[section]
\newtheorem{proposition}[theorem]{Proposition}
\newtheorem{lemma}[theorem]{Lemma}
\newtheorem{corollary}[theorem]{Corollary}
\newtheorem{theoremA}{Theorem}

\newtheorem{theoremB}{Theorem}

\theoremstyle{definition}
\newtheorem{definition}[theorem]{Definition}
\newtheorem{example}[theorem]{Example}
\theoremstyle{remark}
\newtheorem{remark}[theorem]{Remark}
\title[\texorpdfstring{$\delta$}{delta}-Derivations of Finite-Dimensional Simple Lie Algebras]
{\texorpdfstring{$\delta$}{delta}-Derivations of Finite-Dimensional Simple Lie Algebras
with Values in Arbitrary Modules}
\author[Lulu Jiang]{Lulu Jiang}
\address{School of Mathematics, Jilin University, Changchun 130012, Jilin, China}
\email{jianglulu@163.com}

\author[Dong Liu]{Dong Liu}
\address{Department of Mathematics, Huzhou Normal University, Huzhou 313000, Zhejiang, China}
\email{liudong@huznu.edu.cn}

\author{Yufeng Pei}
\address{Department of Mathematics, Huzhou Normal University, Huzhou 313000,
Zhejiang, China}
\email{pei@huznu.edu.cn}
\date{}
\hypersetup{pdfauthor={Lulu Jiang, Dong Liu, Yufeng Pei},
pdftitle={delta-Derivations of Finite-Dimensional Simple Lie Algebras with Values in Arbitrary Modules},
pdfsubject={delta-derivations with values in arbitrary modules}}

\begin{document}
\begin{abstract}
We determine the spaces of \(\delta\)-derivations of finite-dimensional
complex simple Lie algebras with values in arbitrary modules for
\(\delta\ne0,1\). In rank at least two, these derivations vanish
unless \(\delta=\frac{1}{2}\), in which case they are precisely the
module homomorphisms from the adjoint module. For \(\mathfrak{sl}_2\), we
characterize the derivation spaces in terms of Casimir eigenspaces and
determine when they are nonzero. When the coefficient module is simple,
every nonzero derivation module is simple except possibly at \(\delta=2\),
where it has a simple essential socle with a finite-dimensional semisimple
quotient. We describe these modules via translation functors and give
explicit formulas for the derivations.
\end{abstract}
\maketitle

\section{Introduction}

Filippov introduced \(\delta\)-derivations in the adjoint case
\cite{Filippov98} and studied them for prime Lie algebras
\cite{Filippov99}. Derivations with \(\delta=-1\), called
antiderivations, had been considered earlier by Hopkins \cite{Hopkins96}.
Related results on generalized derivations and \(\delta\)-derivations
appear in \cite{LL00,Zus10}.
Zohrabi and Zusmanovich \cite{ZZ21} classified \(\delta\)-derivations
of finite-dimensional semisimple Lie algebras with values in
finite-dimensional modules and suggested studying the case of
infinite-dimensional coefficient modules.

For a complex Lie algebra \(\g\), a \(\g\)-module \(M\), and
\(\delta\in\C\), write \(\Der_\delta(\g,M)\) for the space of
\(\delta\)-derivations \(\g\to M\)
(Definition~\ref{def:delta-derivation}).
We determine these spaces when \(\g\) is finite-dimensional and simple,
with no restriction on \(M\). We assume
\(\delta\in\C\setminus\{0,1\}\): the case \(\delta=0\) is trivial
since \(\g\) is perfect, and \(\delta=1\) gives ordinary derivations.
In rank at least two, we prove the following theorem.

\begin{theoremA}[Theorem~\ref{thm:higher-rank-body} below]\label{thm:higher-rank}
Let \(\g\) be a finite-dimensional complex simple Lie algebra of
rank at least two, let \(M\) be any \(\g\)-module, and let
\(\delta\in\C\setminus\{0,1\}\). Then
\begin{equation}\label{eq:main-arbitrary}
\Der_\delta(\g,M)=
\begin{cases}
\Hom_\g(\g,M),&\delta=\frac{1}{2},\\[3pt]
0,&\delta\ne\frac{1}{2}.
\end{cases}
\end{equation}
The action \eqref{eq:action} on this space is trivial. If \(M\) is simple,
the space is nonzero precisely when \(\delta=\frac{1}{2}\) and
\(M\cong\g\). In that case, fixing a module isomorphism
\(\psi:\g\to M\), one has \(\Der_{\frac{1}{2}}(\g,M)=\C\psi\).
\end{theoremA}

For finite-dimensional \(M\), Theorem~\ref{thm:higher-rank} is contained in
\cite[Main Theorem]{ZZ21}. For arbitrary \(M\), we show that every
\(\delta\)-derivation \(D\) defines the commutative \(2\)-cocycle
\(\Phi_D(x,y)=\frac{1}{2}\bigl(xD(y)+yD(x)\bigr)\)
(Lemma~\ref{lem:comm-from-delta}). Combined with the known vanishing
\(Z^2_{\mathrm{comm}}(\g)=0\) in rank at least two
\cite[Theorem~1.1]{DB09}, \cite[Corollary~2.5]{DZ10}, this yields
Theorem~\ref{thm:higher-rank}.
We also give a short root space proof of this vanishing in
Proposition~\ref{prop:zcomm}, making the proof of the theorem self-contained.

In rank one, \(\g\cong\s\). Let \(V(n)\) denote the simple
\(\s\)-module of highest weight \(n\in\Z_{\ge0}\), and normalize
the Casimir element \(\Omega\) to act on \(V(n)\) by \(n(n+2)\).
For an \(\s\)-module \(X\) and \(a\in\C\), write
\(X^{(a)}=\{v\in X:\Omega v=av\}\), and set
\[
c_\delta=\frac{4(1-\delta)}{\delta^2}.
\]

\begin{theoremB}[Theorem~\ref{thm:reduction} and Corollary~\ref{cor:admissible-parameters} below]\label{thm:main-rank-one}
Let \(\delta\in\C\setminus\{0,1\}\), and let \(M\) be an
\(\s\)-module.
\begin{enumerate}
\item There is a natural isomorphism of \(\s\)-modules
\begin{equation}\label{eq:intro-rank-one}
\Der_\delta(\s,M)
\cong\bigl(V(2)\otimes M^{(c_\delta)}\bigr)^{((1-2\delta)c_\delta)},
\end{equation}
where the tensor product carries the diagonal action.

\item If \(\delta\ne\frac{2}{3}\), then
\(\Der_\delta(\s,M)=0\) if and only if \(M^{(c_\delta)}=0\).
If \(\delta=\frac{2}{3}\), then \(\Der_\delta(\s,M)=0\) if and
only if \(M^{(c_\delta)}\) is a possibly zero direct sum of copies
of \(V(1)\).

\item If \(M\) is simple with Casimir scalar \(\kappa\), then
the parameters \(\delta\in\C\setminus\{0,1\}\) for which
\(\Der_\delta(\s,M)\ne0\) are exactly the roots of
\begin{equation}\label{eq:intro-criterion}
\kappa\delta^2+4\delta-4=0,
\end{equation}
except that \(\delta=\frac{2}{3}\) is omitted when \(M\cong V(1)\).
\end{enumerate}
\end{theoremB}

For \(\s\), we write the \(\delta\)-derivation identity as an
eigenvalue equation for a linear operator on \(\Hom_\C(\s,M)\).
A cubic identity for this operator
(Lemma~\ref{lem:casimir}) and the associated spectral projections yield
Theorem~\ref{thm:main-rank-one}.

For simple \(M\), we use translation functors to determine the
structure of every nonzero derivation module
(Theorem~\ref{thm:translation}). The parameters
\(\delta=\frac{2}{3}\) and \(\delta=2\) correspond to translation
onto and out of the singular block, respectively. At \(\delta=2\),
every nonzero derivation module has a simple essential socle,
and its quotient by the socle is a finite direct sum of copies
of \(V(1)\). Such a module need not be simple, as shown in
Proposition~\ref{prop:weight-two} and Example~\ref{ex:nonsimple}.

We then derive explicit formulas for \(\delta\)-derivations with
values in simple \(\s\)-modules, using Block's classification
\cite{Bavula90,Block81} and the notation of \cite{Maz10}.
These formulas give bases of the derivation spaces for weight
modules and a parametrization by evaluation at \(h\) for
non-weight modules.
For finite-dimensional modules, we recover the results of Zohrabi
and Zusmanovich \cite{ZZ21}. Proposition~\ref{prop:evaluation} and
Theorem~\ref{thm:non-weight} apply in particular to the finite-rank
\(\C[h]\)-torsion-free simple modules described in
\cite{GKM26,GNZ26a,GNZ26b,MX26}.

Section~\ref{sec:rigidity} collects the basic properties of
\(\delta\)-derivations and proves Theorem~\ref{thm:higher-rank}.
Section~\ref{sec:rank-one} proves
Theorem~\ref{thm:main-rank-one} and describes the resulting modules
by translation functors. Section~\ref{sec:explicit} gives explicit formulas for derivations
with values in simple weight and non-weight modules.

\section{The higher rank case}\label{sec:rigidity}

Throughout this paper, all vector spaces, Lie algebras, modules, and
tensor products are over \(\C\), and scalars in operator identities
denote scalar multiples of the identity operator.

Let \(\g\) be a Lie algebra, let \(M\) be a \(\g\)-module,
and let \(\delta\in\C\).

\begin{definition}\label{def:delta-derivation}
A linear map \(D:\g\to M\) is called a \(\delta\)-derivation if
\begin{equation}\label{eq:definition}
D([x,y])=\delta\bigl(xD(y)-yD(x)\bigr)
\qquad(x,y\in\g).
\end{equation}
We denote the space of all \(\delta\)-derivations from \(\g\) to \(M\)
by \(\Der_\delta(\g,M)\).
\end{definition}

The natural \(\g\)-action on \(\Hom_\C(\g,M)\) is given by
\begin{equation}\label{eq:action}
(z\cdot D)(x)=zD(x)-D([z,x]).
\end{equation}

\begin{proposition}\label{prop:basic-properties}
The space \(\Der_\delta(\g,M)\) is a \(\g\)-submodule of
\(\Hom_\C(\g,M)\). Every \(\g\)-module homomorphism
\(F:M\to N\) induces a \(\g\)-module homomorphism
\[
\Der_\delta(\g,M)\longrightarrow\Der_\delta(\g,N),
\qquad D\longmapsto F\circ D.
\]
Moreover, \(\Hom_\g(\g,M)\subseteq\Der_{\frac{1}{2}}(\g,M)\),
and the natural \(\g\)-action on this subspace is trivial.
\end{proposition}
\begin{proof}
Define \(\partial_\delta:\Hom_\C(\g,M)\to
\Hom_\C(\Lambda^2\g,M)\) by
\[
(\partial_\delta D)(x,y)
=D([x,y])-\delta\bigl(xD(y)-yD(x)\bigr).
\]
The \(\g\)-action on \(\Hom_\C(\Lambda^2\g,M)\) is given by
\[
(z\cdot\varphi)(x,y)
=z\varphi(x,y)-\varphi([z,x],y)-\varphi(x,[z,y]).
\]
The Jacobi identity and the module identity imply that
\[
\begin{aligned}
\bigl(\partial_\delta(z\cdot D)\bigr)(x,y)
&=z(\partial_\delta D)(x,y)
 -(\partial_\delta D)([z,x],y)\\
&\qquad -(\partial_\delta D)(x,[z,y])
=(z\cdot\partial_\delta D)(x,y).
\end{aligned}
\]
Thus \(\partial_\delta\) is \(\g\)-equivariant, and its kernel
\(\Der_\delta(\g,M)\) is a \(\g\)-submodule.
For a module homomorphism \(F:M\to N\), one has
\[
\partial_\delta(F\circ D)=F\circ\partial_\delta D,
\qquad z\cdot(F\circ D)=F\circ(z\cdot D),
\]
so composition with \(F\) induces the stated homomorphism.
Finally, for \(D\in\Hom_\g(\g,M)\), equivariance gives
\[
xD(y)-yD(x)=2D([x,y]),\qquad z\cdot D=0.\qedhere
\]
\end{proof}

For \(\delta=0\) and \(\delta=1\), the defining identity \eqref{eq:definition} gives
\[
\Der_0(\g,M)=\Hom_\C(\g/[\g,\g],M),\qquad
\Der_1(\g,M)=Z^1(\g,M),
\]
where \(Z^1(\g,M)\) is the space of Lie algebra
\(1\)-cocycles. Thus \(\Der_0(\g,M)=0\) whenever \(\g\)
is perfect. For ordinary derivations, see \cite{Fuks86,Jacobson62};
for derivations of \(\s\) with values in simple modules, see
\cite{Bavula92}.
Henceforth, we assume \(\delta\notin\{0,1\}\).

\begin{lemma}\label{lem:cyclic-compatibility}
Let \(\g\) be a Lie algebra, let \(M\) be a \(\g\)-module,
and suppose \(\delta\notin\{0,1\}\).
Every \(D\in\Der_\delta(\g,M)\) satisfies
\begin{equation}\label{eq:cyclic-compatibility}
[x,y]D(z)+[y,z]D(x)+[z,x]D(y)=0
\qquad(x,y,z\in\g).
\end{equation}
\end{lemma}
\begin{proof}
Write \(\sum_{\mathrm{cyc}}\) for the cyclic sum over \(x,y,z\).
The \(\delta\)-derivation and module identities give
\begin{equation}\label{eq:cyclic-action}
\sum_{\mathrm{cyc}}zD([x,y])
=\delta\sum_{\mathrm{cyc}}\bigl(zxD(y)-zyD(x)\bigr)
=\delta\sum_{\mathrm{cyc}}[x,y]D(z).
\end{equation}
Applying \(D\) to the Jacobi identity, we obtain
\[
\begin{aligned}
0
&=\sum_{\mathrm{cyc}}D([[x,y],z])\\
&=\delta\sum_{\mathrm{cyc}}
  \bigl([x,y]D(z)-zD([x,y])\bigr)\\
&=\delta(1-\delta)\sum_{\mathrm{cyc}}[x,y]D(z).
\end{aligned}
\]
Since \(\delta\ne0,1\), this establishes \eqref{eq:cyclic-compatibility}.
\end{proof}

We prove Theorem~\ref{thm:higher-rank-body} using the vanishing theorem
for commutative \(2\)-cocycles proved below.

\begin{definition}[{\cite{DZ10}}]
Let \(W\) be a vector space. A symmetric bilinear map
\(\varphi:\g\times \g\to W\) is called a commutative
\(2\)-cocycle if
\[
\varphi([x,y],z)+\varphi([y,z],x)+\varphi([z,x],y)=0
\qquad(x,y,z\in \g).
\]
We write \(Z^2_{\mathrm{comm}}(\g)\) for the space of
\(\C\)-valued commutative \(2\)-cocycles on \(\g\).
\end{definition}

\begin{lemma}\label{lem:comm-from-delta}
Let \(\g\) be a Lie algebra, let \(M\) be a \(\g\)-module,
and suppose \(\delta\notin\{0,1\}\).
For any \(D\in\Der_\delta(\g,M)\), the symmetric bilinear map
\(\Phi_D:\g\times \g\to M\) defined by
\begin{equation}\label{eq:Phi}
\Phi_D(x,y)=\frac{1}{2}\bigl(xD(y)+yD(x)\bigr)
\end{equation}
is a commutative \(2\)-cocycle.
\end{lemma}
\begin{proof}
By \eqref{eq:cyclic-action} and
Lemma~\ref{lem:cyclic-compatibility},
\[
\begin{aligned}
\sum_{\mathrm{cyc}}\Phi_D([x,y],z)
&=\frac{1}{2}\sum_{\mathrm{cyc}}
  \bigl([x,y]D(z)+zD([x,y])\bigr)\\
&=\frac{1+\delta}{2}\sum_{\mathrm{cyc}}[x,y]D(z)=0.
\end{aligned}\qedhere
\]
\end{proof}

\begin{proposition}[{\cite[Corollary~2.5]{DZ10}}]\label{prop:zcomm}
Let \(\g\) be a finite-dimensional complex simple Lie algebra
of rank at least two. Then \(Z^2_{\mathrm{comm}}(\g)=0\).
\end{proposition}
\begin{proof}
The proof in \cite{DZ10} uses the classification of simple
Lie algebras with a subalgebra of codimension one. For the reader's
convenience, we give a direct root space proof in the present setting.
We use the standard root space facts in \cite{Humphreys72}.
Fix a Cartan subalgebra \(\mathfrak h\) with root system \(R\)
and a Weyl-invariant inner product on its real span.
Let \(\varphi\in Z^2_{\mathrm{comm}}(\g)\).
For \(x\in\g_\alpha\), \(y\in\g_\beta\), and
\(h,h'\in\mathfrak h\), the cocycle identity gives
\begin{align}
\alpha(h')\varphi(h,x)&=\alpha(h)\varphi(h',x),\notag\\
(\alpha-\beta)(h)\varphi(x,y)&=-\varphi(h,[x,y]).
\label{eq:zc-root-identity}
\end{align}
Every root \(\gamma\) is a sum of two linearly independent roots:
by irreducibility and the rank assumption, choose \(\alpha\) not proportional
to \(\gamma\) with \((\gamma,\alpha)>0\), and put
\(\beta=\gamma-\alpha\in R\) by the root-string property.
For this pair, \eqref{eq:zc-root-identity} shows that
\(h\mapsto\varphi(h,[x,y])\) is proportional to both
\(\alpha+\beta\) and \(\alpha-\beta\).
These linear forms are independent, so the functional vanishes.
Since \([\g_\alpha,\g_\beta]=\g_\gamma\), we obtain
\begin{equation}\label{eq:zc-Cartan-root}
\varphi(\mathfrak h,\g_\gamma)=0\qquad(\gamma\in R).
\end{equation}

The natural action
\((z\cdot\varphi)(u,v)=-\varphi([z,u],v)-\varphi(u,[z,v])\)
preserves \(Z^2_{\mathrm{comm}}(\g)\). Indeed,
\(z\cdot\varphi\) is symmetric, and the Jacobi identity gives
\[
\begin{aligned}
\sum_{\mathrm{cyc}}(z\cdot\varphi)([u,v],w)
&=-\sum_{\mathrm{cyc}}\bigl(
\varphi([[z,u],v],w)+\varphi([v,w],[z,u])\\
&\hspace{35mm}+\varphi([w,[z,u]],v)\bigr)=0.
\end{aligned}
\]
The sums run over cyclic permutations of \(u,v,w\),
and each expression in parentheses is zero by the cocycle identity.
For arbitrary \(x\in\g_\alpha\) and \(y\in\g_\beta\),
applying \eqref{eq:zc-Cartan-root} to \(x\cdot\varphi\) gives
\[
0=(x\cdot\varphi)(h,y)
=\alpha(h)\varphi(x,y)-\varphi(h,[x,y])
=(2\alpha-\beta)(h)\varphi(x,y),
\]
where the last equality is \eqref{eq:zc-root-identity}.
Since \(R\) is reduced, \(2\alpha-\beta\ne0\), and hence
\(\varphi(\g_\alpha,\g_\beta)=0\) for all roots.
Equation~\eqref{eq:zc-root-identity} also gives
\(\varphi(\mathfrak h,\mathfrak h)=0\), because
\(\mathfrak h\) is spanned by \([\g_\alpha,\g_{-\alpha}]\).
Together with \eqref{eq:zc-Cartan-root}, this proves \(\varphi=0\).
\end{proof}

\begin{theorem}\label{thm:higher-rank-body}
Let \(\g\) be a finite-dimensional complex simple Lie algebra of
rank at least two, let \(M\) be any \(\g\)-module, and let
\(\delta\in\C\setminus\{0,1\}\). Then
\[
\Der_\delta(\g,M)=
\begin{cases}
\Hom_\g(\g,M),&\delta=\frac{1}{2},\\[3pt]
0,&\delta\ne\frac{1}{2}.
\end{cases}
\]
The action \eqref{eq:action} on this space is trivial. If \(M\) is simple,
the space is nonzero precisely when \(\delta=\frac{1}{2}\) and
\(M\cong\g\). In that case, fixing a module isomorphism
\(\psi:\g\to M\), one has \(\Der_{\frac{1}{2}}(\g,M)=\C\psi\).
\end{theorem}
\begin{proof}
Let \(D\in\Der_\delta(\g,M)\). For every \(\ell\in M^*\),
Lemma~\ref{lem:comm-from-delta} and Proposition~\ref{prop:zcomm}
imply that \(\ell\circ\Phi_D=0\). Since the algebraic dual \(M^*\)
separates points, we have \(\Phi_D=0\). Thus \(xD(y)=-yD(x)\), and
\eqref{eq:definition} yields
\begin{equation}\label{eq:higher-rank-action}
xD(y)=\frac{1}{2\delta}D([x,y])\qquad(x,y\in\g).
\end{equation}
When \(\delta=\frac{1}{2}\), equation~\eqref{eq:higher-rank-action}
says that \(D\) is a \(\g\)-module homomorphism. The reverse inclusion
and the triviality of the natural action follow from
Proposition~\ref{prop:basic-properties}.

Suppose now that \(\delta\ne\frac{1}{2}\). Using \eqref{eq:higher-rank-action}
and the Jacobi identity, we obtain
\[
\begin{aligned}
[x,y]D(z)
&=x\bigl(yD(z)\bigr)-y\bigl(xD(z)\bigr)\\
&=\frac{1}{2\delta}\bigl(xD([y,z])-yD([x,z])\bigr)\\
&=\frac{1}{4\delta^2}D\bigl([x,[y,z]]-[y,[x,z]]\bigr)\\
&=\frac{1}{4\delta^2}D([[x,y],z]).
\end{aligned}
\]
On the other hand, \eqref{eq:higher-rank-action} also gives
\([x,y]D(z)=(2\delta)^{-1}D([[x,y],z])\). Hence
\[
(2\delta-1)D([[x,y],z])=0\qquad(x,y,z\in\g).
\]
Since \(\g\) is perfect, \([[\g,\g],\g]=\g\), and hence \(D=0\).

Finally, if \(M\) is simple, any nonzero homomorphism \(\g\to M\)
is an isomorphism because the adjoint module \(\g\) is simple.
For a fixed isomorphism \(\psi:\g\to M\), Schur's lemma gives
\[
\Hom_\g(\g,M)=\psi\circ\End_\g(\g)=\C\psi.\qedhere
\]
\end{proof}

\begin{remark}\label{rem:higher-rank-literature}
For finite-dimensional \(M\), Theorem~\ref{thm:higher-rank-body} is
contained in \cite[Main Theorem]{ZZ21}, including the case
\(\delta=\frac{1}{2}\). For arbitrary \(M\), the vanishing for
\(\delta\notin\{0,1,\frac{1}{2}\}\) also follows from
\cite[Lemma~4.4]{Zus10} and \cite[proof of Corollary~2.5]{DZ10}.
This argument does not cover \(\delta=\frac{1}{2}\), which is
excluded in the cited lemma.

The commutative \(2\)-cocycle argument above also treats
\(\delta=\frac{1}{2}\) for arbitrary \(M\).
Lemma~\ref{lem:comm-from-delta}, combined with the known vanishing
\(Z^2_{\mathrm{comm}}(\g)=0\)
\cite[Theorem~1.1]{DB09}, \cite[Corollary~2.5]{DZ10}, gives
\(\Der_{\frac{1}{2}}(\g,M)=\Hom_\g(\g,M)\).
Proposition~\ref{prop:zcomm} supplies a direct root space proof
of this vanishing, so the proof is self-contained.
\end{remark}

\begin{remark}\label{rem:geometric-interpretation}
For a geometric interpretation, let \(G\) be a connected complex Lie group
with finite-dimensional
Lie algebra \(\g\) and left Maurer--Cartan form \(\theta\).
For a representation \(\rho:\g\to\End_\C M\), set
\(\mathcal E=\mathcal O_G\otimes_\C M\), where \(\mathcal O_G\)
is the sheaf of holomorphic functions on \(G\). Equip \(TG\) with the
left-invariant connection \(\nabla^G_x y=\frac{1}{2}[x,y]\) for
left-invariant vector fields \(x,y\in\g\), and endow \(\mathcal E\) with
the connection \(\nabla^\delta=d+\delta\rho(\theta)\).
Let \(\widehat\nabla^\delta\) be the induced connection on
\(T^*G\otimes\mathcal E\), and let \(d^{\nabla^\delta}\) denote
the covariant exterior derivative. For a linear map \(D:\g\to M\),
define the left-invariant \(\mathcal E\)-valued \(1\)-form
\(\alpha_D=D\circ\theta\). Evaluating on left-invariant vector fields
\(x,y\in\g\) yields
\[
\begin{aligned}
(d^{\nabla^\delta}\alpha_D)(x,y)
&=\delta\bigl(xD(y)-yD(x)\bigr)-D([x,y]),\\
(\widehat\nabla^\delta_x\alpha_D)(y)
&=\delta xD(y)-\frac{1}{2}D([x,y]).
\end{aligned}
\]
Thus \(D\) is a \(\delta\)-derivation precisely when \(\alpha_D\)
is covariantly closed, or equivalently when it satisfies the Codazzi equation
\[
(\widehat\nabla^\delta_x\alpha_D)(y)
=(\widehat\nabla^\delta_y\alpha_D)(x).
\]

At \(\delta=\frac{1}{2}\), the form \(\alpha_D\) is parallel if and only if
\(D\in\Hom_\g(\g,M)\). By Theorem~\ref{thm:higher-rank-body},
every such left-invariant Codazzi form is parallel when \(\g\) is simple of
rank at least two.
\end{remark}

\section{The rank one case}\label{sec:rank-one}

For \(\s\), the derivation identity is an eigenvalue equation
on \(\Hom_\C(\s,M)\). We use this equation to determine when the
derivation space is nonzero and to describe its module structure via translation functors.

\subsection{Casimir identities}\label{sec:casimir}

Fix a basis \(e,h,f\) of \(\s\) with relations
\begin{equation}\label{eq:sl2}
[h,e]=2e,\qquad [h,f]=-2f,\qquad [e,f]=h.
\end{equation}
We use the Casimir element
\[
\Omega=2(ef+fe)+h^2.
\]
The Casimir element used in \cite{Maz10} is \(\Omega+1\).
For the remainder of the paper, set
\begin{equation}\label{eq:parameters}
t=\frac{2}{\delta}-1,\qquad c=t^2-1,\qquad c'=(t-2)^2-1.
\end{equation}
Thus \(\delta=\frac{2}{t+1}\), \(c=c_\delta\), and
\(c'=(1-2\delta)c_\delta\). The values
\(\delta=2,\frac{2}{3}\) correspond to \(t=0,2\), respectively.
The assumption \(\delta\ne0,1\) gives \(t\notin\{-1,1\}\).

For an \(\s\)-module \(X\) and \(a\in\C\), write
\(X^{(a)}=\{v\in X:\Omega v=av\}\).
As above, \(V(n)\) denotes the simple module of highest weight
\(n\in\Z_{\ge0}\).

For an arbitrary \(\s\)-module \(M\), put
\(H=\Hom_\C(\s,M)\). Denote by \(\Omega_M\) the Casimir operator
on \(M\), and also the operator on \(H\) given by
\((\Omega_MD)(x)=\Omega D(x)\). We reserve \(\Omega\) on \(H\)
for the Casimir operator associated with the natural action
\eqref{eq:action}. The linear map \(\Lambda^2\s\to\s\),
\(x\wedge y\mapsto[x,y]\), is an isomorphism: it is surjective
because \([\s,\s]=\s\), and both spaces have dimension three.
Thus the identity
\begin{equation}\label{eq:Sdef}
(\mathrm{S}D)([x,y])=xD(y)-yD(x)
\end{equation}
defines a linear operator \(\mathrm{S}\) on \(H\). By definition,
\(D\) is a \(\delta\)-derivation if and only if
\(\mathrm{S}D=\delta^{-1}D\), or equivalently,
\(2\mathrm{S}D=(t+1)D\). The bracket and the module action are
\(\s\)-equivariant, so \(\mathrm{S}\) is also equivariant.
In the coordinates \(D\mapsto(D(e),D(h),D(f))^{\top}\),
its matrix is
\begin{equation}\label{eq:Smatrix}
\mathrm{S}=\begin{pmatrix}
\frac{h}{2}&-\frac{e}{2}&0\\
-f&0&e\\
0&\frac{f}{2}&-\frac{h}{2}
\end{pmatrix}.
\end{equation}
Each entry is the action of an element of \(U(\s)\), so centrality
of the Casimir element implies \(\Omega_M\mathrm{S}=\mathrm{S}\Omega_M\).
The derivation identity is bilinear and alternating, so it suffices
to check it on \((h,e),(h,f),(e,f)\), which gives
\begin{equation}\label{eq:basis-equations}
\begin{aligned}
(h-t-1)D(e)&=eD(h),\\
(h+t+1)D(f)&=fD(h),\\
(t+1)D(h)&=2eD(f)-2fD(e).
\end{aligned}
\end{equation}

\begin{lemma}\label{lem:casimir}
On \(\Hom_\C(\s,M)\), one has
\begin{equation}\label{eq:operator-identities}
(\mathrm{S}-1)(4\mathrm{S}^2-4\mathrm{S}-\Omega_M)=0,\qquad \Omega=\Omega_M+8-8\mathrm{S}.
\end{equation}
As operators on \(M\), the sums of the diagonal entries in
the matrices of \(\mathrm{S}\) and \(\mathrm{S}^2\) are \(0\) and
\(\frac{\Omega_M}{2}\), respectively.
\end{lemma}
\begin{proof}
To factor \(4\mathrm{S}^2-4\mathrm{S}-\Omega_M\) through \(M\), define
\(\iota:M\to H\) and \(\pi:H\to M\) by
\[
(\iota m)(x)=xm,\qquad
\pi(D)=2fD(e)+hD(h)+2eD(f).
\]
The matrices of \(\iota\) and \(\pi\) are the column
\((e,h,f)^{\top}\) and the row \((2f,h,2e)\).
The module identity gives
\((\mathrm{S}\iota m)([x,y])=(xy-yx)m=[x,y]m\), so \(\mathrm{S}\iota=\iota\).
For the other two compositions, the Lie relations give
\[
\pi \mathrm{S}=(fh-hf,\ ef-fe,\ he-eh)=(2f,h,2e)=\pi,
\qquad \pi\iota=2fe+h^2+2ef=\Omega_M.
\]
Matrix multiplication in \eqref{eq:Smatrix} gives
\[
4\mathrm{S}^2-4\mathrm{S}-\Omega_M=
-\begin{pmatrix}
2ef&eh&2e^2\\
2hf&h^2&2he\\
2f^2&fh&2fe
\end{pmatrix}.
\]
For example, the \((1,1)\)-entry is
\(h^2+2ef-2h-\Omega_M=-2ef\), and the \((1,2)\)-entry
is \(-he+2e=-eh\). Thus
\begin{equation}\label{eq:factor}
4\mathrm{S}^2-4\mathrm{S}-\Omega_M=-\iota\pi.
\end{equation}
Since \((\mathrm{S}-1)\iota=0\), multiplying
\eqref{eq:factor} on the left by \(\mathrm{S}-1\) proves the first identity.

The natural action on \(H\) is represented by
\[
h\mapsto\begin{pmatrix}h-2&0&0\\0&h&0\\0&0&h+2\end{pmatrix},\quad
e\mapsto\begin{pmatrix}e&0&0\\2&e&0\\0&-1&e\end{pmatrix},\quad
f\mapsto\begin{pmatrix}f&1&0\\0&f&-2\\0&0&f\end{pmatrix}.
\]
Substitution into \(\Omega=4fe+h^2+2h\) gives the matrix
\[
\begin{pmatrix}
\Omega_M+8-4h&4e&0\\
8f&\Omega_M+8&-8e\\
0&-4f&\Omega_M+8+4h
\end{pmatrix}
=\Omega_M+8-8\mathrm{S},
\]
which proves the second identity. The diagonal entries of \(\mathrm{S}^2\) are
\[
\frac{h^2}{4}+\frac{ef}{2},\qquad \frac{fe+ef}{2},\qquad \frac{fe}{2}+\frac{h^2}{4},
\]
whose sum is \(\frac{\Omega_M}{2}\). The sum of the diagonal entries
of \(\mathrm{S}\) is zero by \eqref{eq:Smatrix}.
\end{proof}

The invariant form normalized by \((e,f)=1\) and \((h,h)=2\)
identifies \(H\) with \(\s^*\otimes M\cong V(2)\otimes M\).

\begin{remark}\label{rem:characteristic-identity}
For dual bases \((x_i)\) and \((x^i)\) with respect to this form,
the split Casimir satisfies
\[
-2\mathrm{S}=\sum_i\ad(x_i)\otimes x^i.
\]
Indeed, the operator on the right sends \(D\) to the map
\(x\mapsto-\sum_i x^iD([x_i,x])\), so the identity follows
from \eqref{eq:Smatrix}.
Characteristic identities for finite-dimensional modules and for
modules with a central character are treated in
\cite{BrackenGreen71} and \cite{Kostant75}, respectively.
The calculation above applies to arbitrary modules.
\end{remark}

\subsection{Reduction to Casimir eigenspaces}\label{subsec:rank-one-space}

\begin{theorem}\label{thm:reduction}
Let \(\delta\in\C\setminus\{0,1\}\), let \(M\) be an \(\s\)-module,
and put \(c_\delta=\frac{4(1-\delta)}{\delta^2}\).
There is a natural isomorphism of \(\s\)-modules
\begin{equation}\label{eq:rank-one-space}
\Der_\delta(\s,M)
\cong\bigl(V(2)\otimes M^{(c_\delta)}\bigr)^{((1-2\delta)c_\delta)},
\end{equation}
where the tensor product carries the diagonal action.
Moreover:
\begin{enumerate}[nosep]
\item If \(\delta\ne\frac{2}{3}\), then
\(\Der_\delta(\s,M)=0\) if and only if \(M^{(c_\delta)}=0\).

\item If \(\delta=\frac{2}{3}\), then
\(\Der_\delta(\s,M)=0\) if and only if \(M^{(c_\delta)}\) is a possibly zero
direct sum of copies of \(V(1)\).
\end{enumerate}
\end{theorem}
\begin{proof}
Put \(N=M^{(c_\delta)}\). In the notation of \eqref{eq:parameters},
\(c=c_\delta\) and \(c'=(1-2\delta)c_\delta\).
Since the Casimir element is central, \(N\) is a submodule.
For \(D\in\Der_\delta(\s,M)\), use
\(\mathrm{S}D=\frac{t+1}{2}D\) and
\(\Omega_M\mathrm{S}=\mathrm{S}\Omega_M\) in
\eqref{eq:operator-identities} to obtain
\[
\begin{aligned}
0
&=2(\mathrm{S}-1)(4\mathrm{S}^2-4\mathrm{S}-\Omega_M)D\\
&=(t-1)\bigl((t+1)^2-2(t+1)-\Omega_M\bigr)D\\
&=(t-1)(c-\Omega_M)D.
\end{aligned}
\]
Since \(t\ne1\), every value of \(D\) lies in \(N\), so
\(\Der_\delta(\s,M)=\Der_\delta(\s,N)\).
On \(\Hom_\C(\s,N)\), the second identity becomes
\[
\Omega-c'=-8\Bigl(\mathrm{S}-\frac{t+1}{2}\Bigr),
\qquad c+8-c'=4(t+1).
\]
Thus \(\Omega-c'\) and \(\mathrm{S}-\frac{t+1}{2}\) have the same kernel.
The invariant form fixed above gives
\(\Hom_\C(\s,N)\cong\s^*\otimes N\cong V(2)\otimes N\).
Restricting this identification to the \(c'\)-eigenspace of \(\Omega\) gives
\begin{equation}\label{eq:reduction}
\Der_\delta(\s,M)=\Der_\delta(\s,N)
\cong\bigl(V(2)\otimes N\bigr)^{(c')}.
\end{equation}
The basis dual to \(e,h,f\) is \(f,\frac{h}{2},e\), so this map is
\begin{equation}\label{eq:tensor-identification}
D\longmapsto f\otimes D(e)+\frac{1}{2}h\otimes D(h)+e\otimes D(f).
\end{equation}
The inverse sends \(f\otimes u+\frac{1}{2}h\otimes v+e\otimes w\)
to the map with values \(D(e)=u\), \(D(h)=v\), and \(D(f)=w\).
These identifications are \(\s\)-equivariant and commute with
homomorphisms in the second factor, so the isomorphism is natural.

It remains to determine when \(\ker(\mathrm{S}-\frac{t+1}{2})\) is nonzero.
On \(\Hom_\C(\s,N)\), the first identity in
\eqref{eq:operator-identities} factors as
\begin{equation}\label{eq:three-roots}
(\mathrm{S}-1)\Bigl(\mathrm{S}-\frac{t+1}{2}\Bigr)
\Bigl(\mathrm{S}-\frac{1-t}{2}\Bigr)=0.
\end{equation}
Suppose first that \(t\ne0\), and put
\begin{equation}\label{eq:projector}
\mathrm{P}=\frac{(\mathrm{S}-1)(2\mathrm{S}+t-1)}{t(t-1)}.
\end{equation}
The three roots in \eqref{eq:three-roots} are distinct because
\(t\notin\{-1,0,1\}\). The polynomial used in
\eqref{eq:projector} takes the values \(0,1,0\) at
\(1,\frac{t+1}{2},\frac{1-t}{2}\), respectively. Hence
\[
\mathrm{P}^2=\mathrm{P},\qquad
\Bigl(\mathrm{S}-\frac{t+1}{2}\Bigr)\mathrm{P}=0,\qquad
\mathrm{P}D=D\quad
\left(D\in\ker\Bigl(\mathrm{S}-\frac{t+1}{2}\Bigr)\right).
\]
Thus \(\im \mathrm{P}=\ker(\mathrm{S}-\frac{t+1}{2})\).
Write \(\mathrm{P}=(\mathrm{P}_{ij})\) as a matrix of endomorphisms of \(N\).
By Lemma~\ref{lem:casimir},
\begin{equation}\label{eq:partial-trace}
\mathrm{P}_{11}+\mathrm{P}_{22}+\mathrm{P}_{33}
=\frac{c-3(t-1)}{t(t-1)}\,\mathrm{id}_N
=\frac{t-2}{t}\,\mathrm{id}_N.
\end{equation}
If \(N\ne0\) and \(t\notin\{0,2\}\), this sum is nonzero.
Hence \(\mathrm{P}\ne0\) and \(\Der_\delta(\s,N)\ne0\).

Suppose \(t=2\), that is, \(\delta=\frac{2}{3}\), so \(c=3\).
Here \(\mathrm{P}=\mathrm{S}^2-\frac{1}{2}\mathrm{S}-\frac{1}{2}\),
and the derivation space is zero if and only if \(\mathrm{P}=0\).
If \(\mathrm{P}=0\), its \((1,3)\)- and \((3,1)\)-entries,
\(-\frac{e^2}{2}\) and \(-\frac{f^2}{2}\), give \(e^2=f^2=0\) on \(N\).
The identities
\[
0=[f,e^2]=-he-eh,\qquad 0=[e,f^2]=hf+fh
\]
together with \eqref{eq:sl2} give \(he=e\) and \(hf=-f\).
Multiplying \(ef-fe=h\)
on the right by \(e\) and \(f\) then gives
\[
efe=e,\qquad fef=f.
\]
Since \(e^2=f^2=0\), it follows that
\[
h^2=(ef-fe)^2=efef+fefe=ef+fe.
\]
Thus \(\Omega_M=3h^2\), and \(\Omega_M=3\) gives \(h^2=1\).
Consequently,
\[
ef=\frac{1+h}{2},\qquad fe=\frac{1-h}{2},\qquad
N=N_{-1}\oplus N_1,
\]
where \(N_{\pm1}=\{v\in N:hv=\pm v\}\).
The maps \(e:N_{-1}\to N_1\) and
\(f:N_1\to N_{-1}\) are inverse. Choose a basis
\(\{v_i:i\in I\}\) of \(N_{-1}\). Since
\[
fv_i=f^2ev_i=0,\qquad fev_i=v_i,\qquad e^2v_i=0,
\]
each \(\C v_i\oplus\C ev_i\) is a submodule isomorphic to
\(V(1)\), and
\[
N=\bigoplus_{i\in I}(\C v_i\oplus\C ev_i)
\cong\bigoplus_{i\in I}V(1).
\]
Conversely, \(V(2)\otimes V(1)\cong V(3)\oplus V(1)\)
has no Casimir eigenvalue \(-1=c'\). Tensor products and
Casimir eigenspaces commute with direct sums, so
\(\Der_{\frac{2}{3}}(\s,N)=0\) for every direct sum of copies of \(V(1)\).

It remains to consider \(t=0\), that is, \(\delta=2\).
Here \(c=-1\), and \eqref{eq:three-roots} reads
\((\mathrm{S}-1)(\mathrm{S}-\frac{1}{2})^2=0\).
Since \(\frac{1}{2}\) is a repeated root, we replace the projection
by an operator whose image lies in \(\ker(\mathrm{S}-\frac{1}{2})\).
Set
\begin{equation}\label{eq:nilpotent-operator}
\mathrm{J}=-2\mathrm{S}^2+3\mathrm{S}-1
=-2(\mathrm{S}-1)\Bigl(\mathrm{S}-\frac{1}{2}\Bigr).
\end{equation}
The factorization gives
\[
\mathrm{J}^2=4(\mathrm{S}-1)^2\Bigl(\mathrm{S}-\frac{1}{2}\Bigr)^2=0,\qquad
\Bigl(\mathrm{S}-\frac{1}{2}\Bigr)\mathrm{J}
=-2(\mathrm{S}-1)\Bigl(\mathrm{S}-\frac{1}{2}\Bigr)^2=0.
\]
By Lemma~\ref{lem:casimir}, the sum of the diagonal entries
of \(\mathrm{J}\) is
\[
(-c-3)\,\mathrm{id}_N=-2\,\mathrm{id}_N.
\]
If \(N\ne0\), then \(\mathrm{J}\ne0\) and
\[
0\ne\im \mathrm{J}\subseteq
\ker\Bigl(\mathrm{S}-\frac{1}{2}\Bigr)=\Der_2(\s,N).\qedhere
\]
\end{proof}

Since \(U(\s)\) has countable dimension, Dixmier's lemma gives
\(\End_{\s}M=\C\) for every simple \(\s\)-module \(M\);
in particular, \(\Omega\) acts by a scalar \cite{Maz10}.

\begin{corollary}\label{cor:admissible-parameters}
Let \(M\) be a simple \(\s\)-module with Casimir scalar \(\kappa\).
The parameters \(\delta\) for which \(\Der_\delta(\s,M)\ne0\)
are exactly the roots in \(\C\setminus\{0,1\}\) of
\begin{equation}\label{eq:parameter-quadratic}
\kappa\delta^2+4\delta-4=0,
\end{equation}
except that \(\delta=\frac{2}{3}\) is omitted when \(M\cong V(1)\).
There are no such parameters for \(\kappa=0\), and only
\(\delta=2\) for \(\kappa=-1\).
For \(\kappa\ne0\), the two roots, counted with multiplicity,
are interchanged by
\begin{equation}\label{eq:paired-parameters}
\delta\longmapsto\frac{\delta}{\delta-1}.
\end{equation}
\end{corollary}
\begin{proof}
The space \(M^{(c_\delta)}\) is \(M\) when
\(\kappa=c_\delta=\frac{4(1-\delta)}{\delta^2}\), and is zero otherwise.
Theorem~\ref{thm:reduction} gives \eqref{eq:parameter-quadratic}
and the stated exception. For \(\kappa=0\), its only root is the excluded
value \(1\); for \(\kappa=-1\), it has the double root \(2\).
For \(\kappa\ne0\), Vieta's formulas show that the reciprocals of the
roots sum to \(1\), giving \eqref{eq:paired-parameters}.
\end{proof}

\begin{corollary}\label{cor:hom-decomposition}
Let \(\delta\notin\{0,1,2\}\), and let \(M\) be an \(\s\)-module
on which \(\Omega\) acts by \(\frac{4(1-\delta)}{\delta^2}\).
There is a natural direct sum decomposition of \(\s\)-modules
\begin{equation}\label{eq:hom-decomposition}
\Hom_\C(\s,M)
=\Der_1(\s,M)\oplus\Der_\delta(\s,M)
\oplus\Der_{\frac{\delta}{\delta-1}}(\s,M).
\end{equation}
\end{corollary}
\begin{proof}
On \(\Hom_\C(\s,M)\), the first identity in
\eqref{eq:operator-identities} becomes
\[
(\mathrm{S}-1)(\mathrm{S}-\delta^{-1})
\bigl(\mathrm{S}-(1-\delta^{-1})\bigr)=0.
\]
The three roots are distinct under the stated assumptions.
Since \(\ker(\mathrm{S}-\varepsilon^{-1})=\Der_\varepsilon(\s,M)\)
for every \(\varepsilon\ne0\), the eigenspace decomposition gives
\eqref{eq:hom-decomposition}.
The operator \(\mathrm{S}\) is \(\s\)-equivariant and commutes with
maps induced by module homomorphisms, so this
decomposition is natural.
\end{proof}

\subsection{Translation functors}\label{sec:translation}

We use translation functors to determine the structure of the derivation
module for simple coefficient modules.

\begin{definition}\label{def:central-character}
For \(\kappa\in\C\), put \(A_\kappa=U(\s)/(\Omega-\kappa)\).
For an \(\s\)-module \(X\), write
\[
X^{[\kappa]}=\bigcup_{r\ge1}
\{v\in X:(\Omega-\kappa)^r v=0\}.
\]
The full subcategory of \(\s\)-modules satisfying \(X=X^{[\kappa]}\)
is called the generalized central-character block at \(\kappa\).
The \(A_\kappa\)-modules are those annihilated by \(\Omega-\kappa\).
With our normalization, a shifted Harish-Chandra parameter for this
block is a choice of square root of \(\kappa+1\).
The block is singular if this parameter is zero, and regular otherwise.
\end{definition}

The translation functors used below are obtained by tensoring
with \(V(2)\) and taking the appropriate generalized Casimir
summand. The source and target parameters can be chosen as
\(t\) and \(t-2\), whose difference is an extremal weight of \(V(2)\).
The Weyl group changes the sign of a shifted parameter;
on \(\delta\), this is the involution in
\eqref{eq:paired-parameters}.

For an \(A_c\)-module \(N\), equation~\eqref{eq:reduction}
identifies \(\Der_\delta(\s,N)\) with the \(c'\)-eigenspace
in \(V(2)\otimes N\). Thus the derivation functor is a translation
functor whenever this eigenspace equals the corresponding generalized
eigenspace. We use the identification \eqref{eq:tensor-identification}
throughout this subsection.

\begin{theorem}\label{thm:translation}
Let \(\delta\in\C\setminus\{0,1\}\), and put
\[
t=\frac{2}{\delta}-1,\qquad c=t^2-1,\qquad c'=(t-2)^2-1.
\]
\begin{enumerate}
\item If \(\delta\notin\{2,\frac{2}{3}\}\), then
\(\Der_\delta(\s,-):A_c\text{-}\mathrm{Mod}
\to A_{c'}\text{-}\mathrm{Mod}\)
is an equivalence of categories given by a translation functor.
In particular, \(\Der_\delta(\s,M)\) is simple for every simple
\(A_c\)-module \(M\).
\item If \(\delta=\frac{2}{3}\), then \(\Der_{\frac{2}{3}}(\s,-)\) is translation
from \(c=3\) onto \(c'=-1\). It induces a bijection between
the isomorphism classes of simple \(A_3\)-modules other than
\(V(1)\) and those of simple \(A_{-1}\)-modules.
\item If \(\delta=2\) and \(M\) is a simple \(A_{-1}\)-module,
put
\begin{equation}\label{eq:out}
\T(M)=\{v\in V(2)\otimes M:(\Omega-3)^2v=0\}.
\end{equation}
Under the identification \eqref{eq:tensor-identification}, one has
\begin{equation}\label{eq:radical}
\Der_2(\s,M)=\rad\T(M).
\end{equation}
Moreover, \(\soc\T(M)\) is the simple essential socle of
\(\Der_2(\s,M)\), and
\[
\Der_2(\s,M)/\soc\T(M)\cong V(1)^{\oplus m}
\qquad\text{for some }m\in\Z_{\ge0}.
\]
\end{enumerate}
\end{theorem}
\begin{proof}
For (i), assume \(t\notin\{0,1,2\}\) and first work with
finitely generated modules. On \(V(2)\otimes N\), where
\(N\) is an \(A_c\)-module, equations
\eqref{eq:three-roots} and \(\Omega=c+8-8\mathrm{S}\) give
\[
(\Omega-c)(\Omega-c')\bigl(\Omega-((t+2)^2-1)\bigr)=0.
\]
Since
\[
c'-c=4(1-t)\ne0,\qquad
c'-\bigl((t+2)^2-1\bigr)=-8t\ne0,
\]
\(\Omega-c'\) is a simple factor of this annihilating polynomial, so
\[
(V(2)\otimes N)^{[c']}=(V(2)\otimes N)^{(c')}
=\Der_\delta(\s,N).
\]
This summand corresponds to the extremal shift \(t\mapsto t-2\).

If \(t\notin\Z\), both shifted parameters are dominant in
the sense of the translation principle, since neither is a
negative integer. If \(t\in\Z\), then \(t\ge3\) or \(t\le-2\);
in the latter case replace \((t,t-2)\) by \((-t,2-t)\).
In either case the parameters are dominant, have trivial
stabilizers, and differ by an extremal weight \(2\) or \(-2\)
of \(V(2)\). The translation principle
\cite[Theorem~4.1]{BG80}, in the form recalled in
\cite[Lemma~12(i)]{CCM21}, gives inverse equivalences between
the two generalized central-character blocks.

We next check that these equivalences restrict to modules on which
the Casimir acts by a scalar. The equality of the eigenspace and generalized
eigenspace shows that the forward functor takes
\(A_c\)-modules to \(A_{c'}\)-modules. For an \(A_{c'}\)-module \(N\),
the reverse functor selects the generalized \(\mathrm{S}\)-eigenspace
for \(\frac{3-t}{2}\), since
\[
\frac{c'+8-c}{8}=\frac{3-t}{2},\qquad
(\mathrm{S}-1)\Bigl(\mathrm{S}-\frac{t-1}{2}\Bigr)
\Bigl(\mathrm{S}-\frac{3-t}{2}\Bigr)=0.
\]
This root differs from \(1\) and \(\frac{t-1}{2}\) because
\(t\ne1,2\). Hence the generalized \(c\)-eigenspace selected
by the reverse functor is also an eigenspace:
\[
N\longmapsto (V(2)\otimes N)^{[c]}
=\bigl(V(2)\otimes N\bigr)^{(c)}.
\]
Tensoring with the finite-dimensional module \(V(2)\) and taking
direct summands preserve finite generation. The unit and counit
of the block equivalences therefore restrict to finitely generated
\(A_c\)- and \(A_{c'}\)-modules.

For an arbitrary module, write \(N=\bigcup_i N_i\) as the directed
union of its finitely generated submodules. Tensoring with \(V(2)\)
commutes with this union, and each vector in a Casimir eigenspace
already belongs to one of the tensor products. Thus
\[
\bigl(V(2)\otimes\bigcup_iN_i\bigr)^{(\kappa)}
=\bigcup_i(V(2)\otimes N_i)^{(\kappa)}
\qquad(\kappa=c,c').
\]
Both functors therefore commute with these unions. By naturality,
the unit and counit extend to isomorphisms on the unions,
proving (i) for arbitrary modules.

For (ii), \(t=2\) and \((c,c')=(3,-1)\). The factor
\(\mathrm{S}-\frac{3}{2}\) occurs only once in
\eqref{eq:three-roots}, so
\[
\Der_{\frac{2}{3}}(\s,N)
=(V(2)\otimes N)^{(-1)}=(V(2)\otimes N)^{[-1]}
\qquad(N\in A_3\text{-}\mathrm{Mod}).
\]
This is translation from the regular parameter \(2\) onto the
singular parameter \(0\). It sends a simple module to zero or
a simple module and gives a bijection between the isomorphism
classes of simple modules with nonzero image and those of simple
modules in the singular block
\cite[Theorem~26]{CCM21}. For simple \(A_3\)-modules,
Corollary~\ref{cor:admissible-parameters} gives
\[
\Der_{\frac{2}{3}}(\s,M)=0\quad\Longleftrightarrow\quad
M\cong V(1),
\]
which proves (ii).

For (iii), equations \eqref{eq:operator-identities} and
\eqref{eq:nilpotent-operator} give
\begin{equation}\label{eq:out-operators}
\begin{gathered}
\Omega=7-8\mathrm{S},\qquad
\T(M)=\ker\Biggl(\Bigl(\mathrm{S}-\frac{1}{2}\Bigr)^2\Biggr),\\
\Der_2(\s,M)=\Bigl\{D\in\T(M):\mathrm{S}D=\frac{1}{2}D\Bigr\}.
\end{gathered}
\end{equation}
The characteristic identity is
\((\mathrm{S}-1)(\mathrm{S}-\frac{1}{2})^2=0\).
Its two factors are coprime, so
\[
V(2)\otimes M=\ker(\mathrm{S}-1)\oplus\T(M).
\]
Thus \(\T(M)\) is the generalized Casimir summand at \(3\).
The shifts \(\pm2\) from the singular parameter \(0\) give
the same regular target character, whereas the zero shift has
Casimir scalar \(-1\). Hence \(\T\) is translation out of the
singular block.
By \cite[Theorem~33]{CCM21}, \(V(2)\otimes M\), and hence
its direct summand \(\T(M)\), has finite length.
By \cite[Theorem~26(i)--(iii)]{CCM21}, its socle \(L=\soc\T(M)\) is simple and
essential, and
\[
\T(M)/\rad\T(M)\cong L,\qquad [\T(M):L]=2,
\]
where the brackets denote composition multiplicities.
In particular, \(\rad\T(M)\) is its unique maximal submodule.
To identify the derivation space with the radical, we compute
the endomorphism algebra of \(\T(M)\) and determine
\(\mathrm J|_{\T(M)}\).

For \(Y\) in the generalized block at \(3\), put
\[
\Theta(Y)=(V(2)\otimes Y)^{[-1]}.
\]
This is translation back onto the singular block. By
\cite[Lemma~12(ii) and Theorem~26(i)]{CCM21},
\[
\Theta\T(M)\cong M^{\oplus2},\qquad \Theta L\cong M.
\]
Biadjointness of the translation functors and \(\End_{\s}M=\C\) give
\[
\begin{aligned}
\dim\End_{\s}\T(M)
&=\dim\Hom_{\s}(M,\Theta\T(M))\\
&=\dim\Hom_{\s}(M,M^{\oplus2})=2.
\end{aligned}
\]
Let \(i:L\hookrightarrow\T(M)\) be the socle inclusion and choose
\(q:\T(M)\twoheadrightarrow L\) with \(\ker q=\rad\T(M)\).
If \(qi\ne0\), it is an automorphism of the simple module
\(L\). Replacing \(i\) by \(i(qi)^{-1}\) would give
\(qi=\mathrm{id}_L\), hence
\(\T(M)=i(L)\oplus\ker q\). Essentiality of \(i(L)\)
would force \(\ker q=0\), contrary to \([\T(M):L]=2\).
Thus \(qi=0\). Since \(q\) is surjective and \(i\) is
injective, \(iq\ne0\), and
\[
(iq)^2=i(qi)q=0,\qquad
\End_{\s}\T(M)=\C\,\mathrm{id}\oplus\C\,iq.
\]

By the proof of
Theorem~\ref{thm:reduction}, \(\mathrm{J}\ne0\) on \(V(2)\otimes M\).
By \eqref{eq:nilpotent-operator}, \(\mathrm{J}\) vanishes on
\(\ker(\mathrm{S}-1)\). On \(\T(M)\), the same factorization gives
\[
\begin{aligned}
\mathrm{J}
&=-2\Biggl[\Bigl(\mathrm{S}-\frac{1}{2}\Bigr)-\frac{1}{2}\Biggr]
\Bigl(\mathrm{S}-\frac{1}{2}\Bigr)\\
&=-2\Bigl(\mathrm{S}-\frac{1}{2}\Bigr)^2
  +\Bigl(\mathrm{S}-\frac{1}{2}\Bigr)
 =\mathrm{S}-\frac{1}{2}.
\end{aligned}
\]
Thus \(\mathrm{J}|_{\T(M)}\ne0\) and \(\mathrm{J}^2=0\).
Write \(\mathrm{J}|_{\T(M)}=\alpha\,\mathrm{id}+\beta iq\).
Since \(qi=0\),
\[
q\mathrm{J}=\alpha q,\qquad
0=q\mathrm{J}^2=\alpha^2q.
\]
Surjectivity of \(q\) gives \(\alpha=0\), and nonvanishing of
\(\mathrm{J}|_{\T(M)}\) gives \(\beta\ne0\). Hence
\[
\im\mathrm{J}=i(L)=\soc\T(M),\qquad
\ker\bigl(\mathrm{J}|_{\T(M)}\bigr)=\ker q=\rad\T(M).
\]
On \(\T(M)\), one has
\(\mathrm{J}=\mathrm{S}-\frac{1}{2}=-\frac{1}{8}(\Omega-3)\).
Together with \eqref{eq:out-operators}, this proves
\eqref{eq:radical} and gives the exact sequence
\begin{equation}\label{eq:translated-exact}
0\longrightarrow\Der_2(\s,M)\longrightarrow\T(M)
\xrightarrow{\ \Omega-3\ }\soc\T(M)\longrightarrow0.
\end{equation}
The first map is the inclusion under the identification
\eqref{eq:tensor-identification}.
Since \(\mathrm{J}^2=0\), the simple socle of \(\T(M)\) lies in
\(\Der_2(\s,M)\). Being essential in \(\T(M)\), it is also
essential in \(\Der_2(\s,M)\).

Finally, put \(Q=\Der_2(\s,M)/L\). The functor \(\Theta\) is exact,
so applying it to \eqref{eq:translated-exact} gives
\[
0\longrightarrow\Theta\Der_2(\s,M)
\longrightarrow M^{\oplus2}\longrightarrow M\longrightarrow0.
\]
Hence \(\Theta\Der_2(\s,M)\cong M\). The inclusion
\(L\hookrightarrow\Der_2(\s,M)\) therefore induces an isomorphism
\(\Theta L\xrightarrow{\sim}\Theta\Der_2(\s,M)\).
Exactness now gives \(\Theta Q=0\).
Since \(\Omega\) acts by \(3\) on \(Q\), part~\textup{(ii)} identifies
\(\Theta Q\) with \(\Der_{\frac{2}{3}}(\s,Q)\).
Theorem~\ref{thm:reduction}\textup{(ii)} implies that \(Q\) is a direct
sum of copies of \(V(1)\). As \(Q\) is a subquotient of the
finite-length module \(\T(M)\), only finitely many copies occur.
\end{proof}

\section{Derivations with values in simple modules}\label{sec:explicit}

We give explicit formulas for \(\delta\)-derivations with values
in simple \(\s\)-modules. Recall that \(t=2/\delta-1\).
By Corollary~\ref{cor:admissible-parameters}, we need only consider
modules with Casimir scalar \(t^2-1\).

\subsection{Weight modules}\label{subsec:weight}

The Verma module \(M(\lambda)\) has a basis \((v_n)_{n\ge0}\) on which
\begin{equation}\label{eq:verma}
hv_n=(\lambda-2n)v_n,\quad fv_n=v_{n+1},\quad
ev_n=n(\lambda-n+1)v_{n-1}.
\end{equation}
It is simple exactly when \(\lambda\notin\Z_{\ge0}\).
Let \(\theta\) be the involution
\[
\theta(e)=f,\qquad \theta(f)=e,\qquad \theta(h)=-h,
\]
and put \(M'(\mu)=M(-\mu)^\theta\), where
\(x\cdot_\theta m=\theta(x)m\). The finite-dimensional module \(V(n)\) has basis
\(v_0,\ldots,v_n\) and action \eqref{eq:verma} with \(\lambda=n\).
A basis vector outside the specified index range is understood to be zero.

For \(\xi\in\C/2\Z\) and \(\tau\in\C\), the dense module \(V(\xi,\tau)\)
has basis \((v_\nu)_{\nu\in\xi}\) with action
\begin{equation}\label{eq:dense}
hv_\nu=\nu v_\nu,\qquad fv_\nu=v_{\nu-2},\qquad
ev_\nu=\frac{\tau-(\nu+1)^2}{4}v_{\nu+2}
\quad(\nu\in\xi).
\end{equation}
Its Casimir scalar is \(\tau-1\), and it is simple if and only if
\(-1\pm\sqrt\tau\notin\xi\). For these simplicity criteria,
see \cite{Maz10}.

For a weight module \(M\), write
\[
\Hom_\C(\s,M)_\mu
=\{D:D(e)\in M_{\mu+2},\ D(h)\in M_\mu,\ D(f)\in M_{\mu-2}\}.
\]
Since \(\s\) is finite-dimensional, each map has only finitely many
nonzero homogeneous components. Equation~\eqref{eq:basis-equations}
holds separately in each degree, so
\[
\Der_\delta(\s,M)=\bigoplus_{\mu\in\C}
\bigl(\Der_\delta(\s,M)\cap\Hom_\C(\s,M)_\mu\bigr).
\]
To obtain bases from highest weight derivations, we record the powers
of the lowering operator. Its action is
\begin{equation}\label{eq:cochain-lowering}
f\cdot(D(e),D(h),D(f))=
\bigl(fD(e)+D(h),\ fD(h)-2D(f),\ fD(f)\bigr).
\end{equation}
The off-diagonal part of this action has scalar entries, its cube
is zero, and it commutes with the componentwise action of \(f\).
The binomial formula therefore gives, for \(k\in\Z_{\ge0}\),
\begin{equation}\label{eq:cochain-lowering-powers}
\begin{aligned}
(f^k\cdot D)(e)&=f^kD(e)+kf^{k-1}D(h)-k(k-1)f^{k-2}D(f),\\
(f^k\cdot D)(h)&=f^kD(h)-2kf^{k-1}D(f),\\
(f^k\cdot D)(f)&=f^kD(f).
\end{aligned}
\end{equation}
Here the powers on the right act on \(M\), and terms with
negative powers are omitted.

For finite-dimensional simple modules, Theorem~\ref{thm:reduction}
and the Clebsch--Gordan decomposition recover the result of
Zohrabi and Zusmanovich below. The basis formulas follow by applying
\eqref{eq:cochain-lowering-powers} to a highest weight vector.

\begin{proposition}[{\cite[Lemma~9(ii),(iii)]{ZZ21}}]\label{prop:finite}
The nonzero derivation spaces with values in \(V(n)\) have the following bases.
\begin{enumerate}
\item If \(n\ge1\) and \(\delta=-\frac{2}{n}\), a basis is
\(D_0,\ldots,D_{n+2}\), where
\begin{equation}\label{eq:lower-basis}
D_k(e)=-k(k-1)v_{k-2},\qquad
D_k(h)=-2kv_{k-1},\qquad D_k(f)=v_k.
\end{equation}
\item If \(n\ge2\) and \(\delta=\frac{2}{n+2}\), a basis is
\(D_0,\ldots,D_{n-2}\), where
\begin{equation}\label{eq:upper-basis}
\begin{aligned}
D_k(e)&=-(n-k)(n-k-1)v_k,\\
D_k(h)&=2(n-k-1)v_{k+1},\qquad D_k(f)=v_{k+2}.
\end{aligned}
\end{equation}
\end{enumerate}
All other spaces vanish, including those with values in
\(V(0)\). Under the natural action,
\begin{equation}\label{eq:finite-table}
\begin{aligned}
\Der_{-\frac{2}{n}}(\s,V(n))&\cong V(n+2) &&(n\ge1),\\
\Der_{\frac{2}{n+2}}(\s,V(n))&\cong V(n-2) &&(n\ge2).
\end{aligned}
\end{equation}
\end{proposition}
\begin{proof}
The Casimir condition gives
\(t^2-1=n(n+2)\), hence \(t=\pm(n+1)\).
For \(n\ge2\),
\[
V(2)\otimes V(n)\cong V(n+2)\oplus V(n)\oplus V(n-2),
\]
with respective \(\mathrm{S}\)-eigenvalues
\(-\frac{n}{2},1,\frac{n+2}{2}\) by \eqref{eq:operator-identities}.
For \(n=1\), the last summand is absent. Taking the
\(\delta^{-1}\)-eigenspace gives \eqref{eq:finite-table}.
For \(V(0)\), the derivation space is zero because \(\s\) is perfect.

In the first family, \(D_0(e)=D_0(h)=0\) and \(D_0(f)=v_0\).
Since \(t=-n-1\), equations \eqref{eq:basis-equations} hold by
\((h+t+1)v_0=0\) and \(ev_0=0\).
In the second family, \(t=n+1\), and \eqref{eq:upper-basis} at
\(k=0\) gives
\[
\begin{aligned}
(h-t-1)D_0(e)&=2n(n-1)v_0=eD_0(h),\\
(h+t+1)D_0(f)&=2(n-1)v_2=fD_0(h),\\
eD_0(f)-fD_0(e)&=(n+2)(n-1)v_1
=\frac{t+1}{2}D_0(h).
\end{aligned}
\]
Thus \(D_0\) is a derivation in both families.
Formula~\eqref{eq:cochain-lowering-powers} gives \(D_k=f^k\cdot D_0\).
The first coordinate in \eqref{eq:upper-basis} follows from
\[
-n(n-1)+2k(n-1)-k(k-1)=-(n-k)(n-k-1).
\]
The weights are \(n+2-2k\) and \(n-2-2k\), respectively.
In the first family, nonvanishing follows from
\[
D_k(f)=v_k\ (0\le k\le n),\qquad
D_{n+1}(h)=-2(n+1)v_n,\qquad
D_{n+2}(e)=-(n+2)(n+1)v_n.
\]
In the second family, \(D_k(f)=v_{k+2}\ne0\) for every index
in the stated range.
The displayed maps are thus nonzero and have distinct weights.
By \eqref{eq:finite-table}, they form bases of the two derivation spaces.
\end{proof}

\begin{corollary}[{\cite[Main Theorem]{ZZ21}}]\label{cor:finite-coefficients}
Let \(\g\) be a finite-dimensional complex simple Lie algebra,
let \(M=\bigoplus_i M_i\) be a finite-dimensional \(\g\)-module
with simple summands, and let \(\delta\notin\{0,1\}\). Then
\begin{equation}\label{eq:finite-coefficients}
\Der_\delta(\g,M)=\bigoplus_i\Der_\delta(\g,M_i),
\qquad
\Der_{\frac{1}{2}}(\g,M)=\Hom_\g(\g,M).
\end{equation}
If \(\delta\ne\frac{1}{2}\), the space is nonzero precisely when
\(\g\cong\s\) and \(M\) contains a simple summand \(V(n)\) with
either \(\delta=-\frac{2}{n}\), \(n\ge1\), or
\(\delta=\frac{2}{n+2}\), \(n\ge3\).
Every derivation is then a sum of maps into these summands, with
bases given in Proposition~\ref{prop:finite}.
\end{corollary}
\begin{proof}
Projection onto the simple summands gives the first equality.
By Theorem~\ref{thm:higher-rank-body} and Proposition~\ref{prop:finite},
only adjoint summands contribute at \(\delta=\frac{1}{2}\), and
on each such summand the derivations are the module homomorphisms
from \(\g\). The same results give the assertion for
\(\delta\ne\frac{1}{2}\).
\end{proof}

For infinite-dimensional simple weight modules and \(\delta\ne2\),
Theorem~\ref{thm:translation} shows that every nonzero derivation
module is simple.

\begin{proposition}\label{prop:weight-regular}
Suppose \(\delta\ne2\).
\begin{enumerate}
\item For derivations with values in infinite-dimensional simple
highest weight modules, the nonzero spaces are
\begin{equation}\label{eq:verma-table}
\begin{aligned}
\Der_\delta(\s,M(t-1))&\cong M(t-3),
&&t-1\notin\Z_{\ge0},\\
\Der_\delta(\s,M(-t-1))&\cong M(1-t),
&&-t-1\notin\Z_{\ge0}.
\end{aligned}
\end{equation}
A basis \((D_k)_{k\ge0}\) is given by \eqref{eq:upper-basis} with
\(n=t-1\) for \(M(t-1)\), and by \eqref{eq:lower-basis} for
\(M(-t-1)\).
\item For derivations with values in infinite-dimensional simple
lowest weight modules, the nonzero spaces are
\begin{equation}\label{eq:lowest-table}
\begin{aligned}
\Der_\delta(\s,M'(1-t))&\cong M'(3-t),
&&1-t\notin\Z_{\le0},\\
\Der_\delta(\s,M'(t+1))&\cong M'(t-1),
&&t+1\notin\Z_{\le0}.
\end{aligned}
\end{equation}
For each \(M\) in \textup{(i)}, a basis is
\((D_k\circ\theta:\s\to M^\theta)_{k\ge0}\).
\item The simple dense modules admitting nonzero derivations are
\(V(\xi,t^2)\), with \(t-1,-t-1\notin\xi\). Under the natural action,
\begin{equation}\label{eq:dense-answer}
\Der_\delta\bigl(\s,V(\xi,t^2)\bigr)
\cong V(\xi,(t-2)^2).
\end{equation}
A basis \((D_\mu)_{\mu\in\xi}\) of the derivation space is given by
\begin{equation}\label{eq:dense-general-D}
\begin{aligned}
D_\mu(e)&=-\frac{(\mu+t+1)(\mu+t-1)}{4}v_{\mu+2},\\
D_\mu(h)&=(\mu+t-1)v_\mu,\qquad
D_\mu(f)=v_{\mu-2}.
\end{aligned}
\end{equation}
\end{enumerate}
\end{proposition}
\begin{proof}
For \(M(\lambda)\), the Casimir equation
\(\lambda(\lambda+2)=t^2-1\) gives \(\lambda=t-1\) or
\(\lambda=-t-1\). For \(\lambda=t-1\), set
\[
D_0(e)=-\lambda(\lambda-1)v_0,\qquad
D_0(h)=2(\lambda-1)v_1,\qquad D_0(f)=v_2.
\]
For \(\lambda=-t-1\), set
\[
D_0(e)=D_0(h)=0,\qquad D_0(f)=v_0.
\]
The verification of \eqref{eq:basis-equations} in the proof of
Proposition~\ref{prop:finite} uses only \eqref{eq:verma} and remains
valid with \(n\) replaced by \(\lambda\).
The two maps have weights \(\lambda-2\) and \(\lambda+2\), respectively,
and in both cases
\[
e\cdot D_0=
\bigl(eD_0(e),\ eD_0(h)+2D_0(e),\ eD_0(f)-D_0(h)\bigr)=0.
\]
Theorem~\ref{thm:translation}(i),(ii) implies that the derivation
space is simple, so it is generated by \(D_0\).
The corresponding Verma modules have highest weights \(\lambda-2\)
and \(\lambda+2\), respectively, and are simple. Indeed,
\(\lambda\notin\Z_{\ge0}\) implies
\(\lambda-2\notin\Z_{\ge0}\), whereas
\(\lambda+2\in\Z_{\ge0}\) would force \(\lambda=-1,-2\),
which give the excluded values \(\delta=2,1\).
The homomorphism from the corresponding Verma module that sends
its highest weight vector to \(D_0\) is therefore an isomorphism.
Formula~\eqref{eq:cochain-lowering-powers} gives the stated basis
\(D_k=f^k\cdot D_0\).

Since \(\theta\) is a Lie involution, twisting gives
\[
\Der_\delta(\s,M)^\theta\xrightarrow{\ \sim\ }
\Der_\delta(\s,M^\theta),\qquad D\longmapsto D\circ\theta,
\]
with \(z\cdot(D\circ\theta)=(\theta(z)\cdot D)\circ\theta\).
Applying this isomorphism to the two highest weight modules proves
\eqref{eq:lowest-table} and the corresponding basis formulas.

A homogeneous map with values in a dense module has the form
\[
D(e)=\alpha v_{\mu+2},\qquad D(h)=\beta v_\mu,\qquad
D(f)=\gamma v_{\mu-2}.
\]
The second and third equations of \eqref{eq:basis-equations} give
\[
\begin{aligned}
\beta&=(\mu+t-1)\gamma,\\
\alpha&=\frac{t^2-(\mu-1)^2}{4}\gamma-\frac{t+1}{2}\beta
=-\frac{(\mu+t+1)(\mu+t-1)}{4}\gamma.
\end{aligned}
\]
Substitution into the first equation gives
\[
\begin{aligned}
(\mu-t+1)\alpha
&=-\frac{(\mu-t+1)(\mu+t+1)(\mu+t-1)}{4}\gamma\\
&=\frac{t^2-(\mu+1)^2}{4}\beta.
\end{aligned}
\]
Hence the weight-\(\mu\) solution space is exactly \(\C D_\mu\).
Evaluation at \(f\) is injective on each homogeneous component of the
derivation space, and
\[
(e\cdot D_\mu)(f)
=\left(\frac{t^2-(\mu-1)^2}{4}-(\mu+t-1)\right)v_\mu
=\frac{(t-2)^2-(\mu+1)^2}{4}v_\mu.
\]
Together with the \(h\)- and \(f\)-actions, this gives
\[
h\cdot D_\mu=\mu D_\mu,\qquad f\cdot D_\mu=D_{\mu-2},\qquad
e\cdot D_\mu=\frac{(t-2)^2-(\mu+1)^2}{4}D_{\mu+2}.
\]
This proves \eqref{eq:dense-answer}. The target module satisfies
the same simplicity criterion, since \(t-3\equiv t-1\) and
\(1-t\equiv-t-1\pmod{2\Z}\).
\end{proof}

\begin{remark}\label{rem:half-specialization}
For infinite-dimensional \(\s\)-modules, the equality
\(\Der_{\frac{1}{2}}(\s,M)=\Hom_{\s}(\s,M)\) in
Corollary~\ref{cor:finite-coefficients} may fail. For example,
Proposition~\ref{prop:weight-regular} and the absence of weight \(2\)
in \(M(-4)\) give
\[
\Der_{\frac{1}{2}}(\s,M(-4))\cong M(-2),\qquad
\Hom_{\s}(\s,M(-4))=0.
\]
\end{remark}

At \(\delta=2\), every nonzero derivation module has a simple
essential socle by Theorem~\ref{thm:translation}(iii).
The highest and lowest weight cases below show that the derivation
module need not be simple.

\begin{proposition}\label{prop:weight-two}
The nonzero \(2\)-derivation spaces with values in simple weight modules are
\begin{equation}\label{eq:weight-two}
\begin{aligned}
\Der_2(\s,M(-1))&\cong M(1),\\
\Der_2(\s,M'(1))&\cong M'(-1),\\
\Der_2(\s,V(\xi,0))&\cong V(\xi,4),
\qquad \xi\ne1+2\Z.
\end{aligned}
\end{equation}
Bases for these spaces are given, respectively, by \eqref{eq:lower-basis}
with \(k\ge0\), by composing those maps with \(\theta\), and by
\eqref{eq:dense-general-D} at \(t=0\).
The first two derivation modules have simple socles \(M(-3)\) and
\(M'(3)\), respectively. In both cases, the quotient by the socle
is \(V(1)\). The derivation module in the dense case is simple.
\end{proposition}
\begin{proof}
By Corollary~\ref{cor:admissible-parameters}, only simple weight
modules with Casimir scalar \(-1\) can occur; these are the modules
listed in \eqref{eq:weight-two}. Every nonzero homogeneous map
\(\s\to M(-1)\) has weight \(1-2k\) for some \(k\ge0\), and has the form
\[
D(e)=\alpha v_{k-2},\qquad D(h)=\beta v_{k-1},\qquad
D(f)=\gamma v_k.
\]
For \(k\ge2\), the last two equations of \eqref{eq:basis-equations} and
\(ev_k=-k^2v_{k-1}\) give
\[
\beta=-2k\gamma,\qquad
\alpha=-k^2\gamma-\frac{1}{2}\beta=-k(k-1)\gamma.
\]
The first equation is then automatically satisfied, since
\[
(2-2k)\alpha=2k(k-1)^2\gamma=-(k-1)^2\beta.
\]
For \(k=0\), only
\(D(f)\) can be nonzero; for \(k=1\),
\(D(e)=0\) and \(\beta=-2\gamma\).
Thus each homogeneous solution is \(\gamma D_k\).
Since
\[
(e\cdot D_k)(f)=eD_k(f)-D_k(h)
=(-k^2+2k)v_{k-1},
\]
the action is
\[
h\cdot D_k=(1-2k)D_k,\qquad f\cdot D_k=D_{k+1},\qquad
e\cdot D_k=k(2-k)D_{k-1}.
\]
Consequently the space is isomorphic to \(M(1)\), with
\[
\soc\Der_2(\s,M(-1))
=\bigoplus_{k\ge2}\C D_k\cong M(-3),\qquad
M(1)/M(-3)\cong V(1).
\]
Twisting by \(\theta\) gives the lowest weight case. The calculation
for dense modules remains valid at \(t=0\). Its target \(V(\xi,4)\)
is simple because \(1\) and \(-3\), the two weights excluded by
the simplicity criterion, lie in \(1+2\Z\).
\end{proof}

\subsection{Non-weight modules}\label{subsec:non-weight}

Put \(R=\C[h]\), \(K=\C(h)\), and
\(\rank_R M=\dim_K(K\otimes_R M)\).
Every simple non-weight \(\s\)-module is \(R\)-torsion-free
\cite[Theorem~6.1.3]{Maz10}. We first treat arbitrary
\(R\)-torsion-free modules with scalar Casimir action; finite
generation over \(R\) is not required.

For an \(R\)-torsion-free \(\s\)-module \(M\), the action extends
to \(K\otimes_R M\) by
\begin{equation}\label{eq:localized-action}
e(r(h)m)=r(h-2)em,\qquad
f(r(h)m)=r(h+2)fm,\qquad h(r(h)m)=hr(h)m.
\end{equation}
Here \(r(h)\in K\) and \(m\in M\).
We identify \(M\) with its image in \(K\otimes_R M\).
Rational factors in the formulas below act on the left in this
localization. We specify the conditions under which the resulting
vectors belong to \(M\).

\begin{proposition}\label{prop:evaluation}
Let \(M\ne0\) be an \(R\)-torsion-free \(\s\)-module on
which \(\Omega\) acts by a scalar \(\kappa\), and put \(t=2/\delta-1\).
Then \(\Der_\delta(\s,M)\ne0\) if and only if \(\kappa=t^2-1\).
When \(\kappa=t^2-1\), put
\begin{equation}\label{eq:lattice}
\Lambda=\{u\in M:eu\in(h-t-1)M,\ fu\in(h+t+1)M\}.
\end{equation}
Evaluation \(D\mapsto D(h)\) is a vector space isomorphism
\(\Der_\delta(\s,M)\to\Lambda\), whose inverse is
\begin{equation}\label{eq:Du}
D_u(e)=(h-t-1)^{-1}eu,\qquad D_u(h)=u,\qquad
D_u(f)=(h+t+1)^{-1}fu.
\end{equation}
Moreover, \(\rank_R\Lambda=\rank_R M\).
The natural action \eqref{eq:action} induces the following action on \(\Lambda\):
\begin{equation}\label{eq:transport}
h*u=hu,\qquad e*u=\frac{h-t+1}{h-t-1}eu,\qquad
f*u=\frac{h+t-1}{h+t+1}fu.
\end{equation}
For this action, \(\Omega\) acts by \((t-2)^2-1\).
\end{proposition}
\begin{proof}
Let \(D\in\Der_\delta(\s,M)\) and put \(u=D(h)\).
The first two equations of \eqref{eq:basis-equations} give
\[
(h-t-1)D(e)=eu,\qquad (h+t+1)D(f)=fu.
\]
Since \(M\) is torsion-free, these equations determine
\(D(e)\) and \(D(f)\) uniquely, as in \eqref{eq:Du}.
Hence evaluation at \(h\) is injective. To obtain the restriction
on \(\kappa\), use the identities
\[
fe=\frac{\kappa-h^2-2h}{4},\qquad
ef=\frac{\kappa-h^2+2h}{4},
\]
and the shift relations
\[
eD(f)=(h+t-1)^{-1}efu,\qquad
fD(e)=(h-t+1)^{-1}feu.
\]
Consequently,
\[
\begin{aligned}
2\bigl(eD(f)-fD(e)\bigr)
&=\left(\frac{\kappa-h^2+2h}{2(h+t-1)}
 -\frac{\kappa-h^2-2h}{2(h-t+1)}\right)u\\
&=\frac{(t+1)h^2-(t-1)\kappa}{(h-t+1)(h+t-1)}u.
\end{aligned}
\]
Subtracting \((t+1)u\) and using the third equation of
\eqref{eq:basis-equations}, we obtain
\begin{equation}\label{eq:non-weight-compatibility}
2\bigl(eD(f)-fD(e)\bigr)-(t+1)u
=\frac{(t-1)(t^2-1-\kappa)}{(h-t+1)(h+t-1)}u=0.
\end{equation}
If \(D\ne0\), injectivity of evaluation gives \(u\ne0\),
so torsion-freeness and \(t\ne1\) imply \(\kappa=t^2-1\).
Now suppose \(\kappa=t^2-1\). The first two equations show that
evaluation takes values in \(\Lambda\). Conversely, for \(u\in\Lambda\),
formula~\eqref{eq:Du} defines a map \(D_u:\s\to M\) satisfying
the first two equations of \eqref{eq:basis-equations}.
The calculation in \eqref{eq:non-weight-compatibility} shows that
it also satisfies the third. Thus evaluation is an isomorphism
onto \(\Lambda\).

The shift relations make \(\Lambda\) an \(R\)-submodule.
For \(m\in M\), they give
\[
\begin{aligned}
e\bigl((h-t+1)(h+t-1)m\bigr)&=(h-t-1)(h+t-3)em,\\
f\bigl((h-t+1)(h+t-1)m\bigr)&=(h-t+3)(h+t+1)fm.
\end{aligned}
\]
Thus
\begin{equation}\label{eq:full-rank}
(h-t+1)(h+t-1)M\subseteq\Lambda\subseteq M.
\end{equation}
Since \(M\ne0\) is torsion-free, the leftmost submodule in
\eqref{eq:full-rank} is nonzero. This proves \(\Lambda\ne0\) and
completes the nonvanishing criterion. Localizing
\eqref{eq:full-rank} gives
\(K\otimes_R\Lambda=K\otimes_R M\), hence the rank equality.

We now compute the transported action. The identities
\[
(e\cdot D_u)(h)=eu+2D_u(e),\qquad
(f\cdot D_u)(h)=fu-2D_u(f),\qquad
(h\cdot D_u)(h)=hu
\]
give \eqref{eq:transport}. Since
\(fe=-\frac{1}{4}(h-t+1)(h+t+1)\) on \(M\), we obtain
\[
f*(e*u)
=\frac{(h+t-1)(h-t+3)}{(h+t+1)(h-t+1)}feu.
\]
It follows that
\begin{equation}\label{eq:non-weight-target-casimir}
f*(e*u)=-\frac{(h+t-1)(h-t+3)}{4}u,\qquad
\Omega*u=((t-2)^2-1)u,
\end{equation}
where the second identity follows from
\(\Omega*u=h^2u+2hu+4f*(e*u)\).
\end{proof}

\begin{theorem}\label{thm:non-weight}
Let \(t=2/\delta-1\), and let \(M\) be a simple non-weight
\(\s\)-module with Casimir scalar \(t^2-1\).
Let \(\Lambda\) be given by \eqref{eq:lattice}, with the action
\eqref{eq:transport}.
\begin{enumerate}
\item If \(\delta\ne2\), then \(\Der_\delta(\s,M)\cong\Lambda\)
is a simple non-weight module.

\item If \(\delta=2\), put
\begin{equation}\label{eq:Lambda0}
\Lambda_0=(h^2-1)M+(h-1)fM+(h+1)eM.
\end{equation}
Then \(\Lambda_0\) is the simple essential socle of \(\Lambda\),
and there is an exact sequence
\begin{equation}\label{eq:non-weight-extension}
0\longrightarrow\Lambda_0\longrightarrow\Lambda
\longrightarrow V(1)^{\oplus m}\longrightarrow0,
\qquad m\in\Z_{\ge0}.
\end{equation}
The maps
\begin{equation}\label{eq:fibers}
\bar e:M/(h+1)M\longrightarrow M/(h-1)M,\qquad
\bar f:M/(h-1)M\longrightarrow M/(h+1)M
\end{equation}
induced by the original actions satisfy
\(\bar e\bar f=\bar f\bar e=0\), and
\begin{equation}\label{eq:multiplicity}
m=\dim_\C\frac{\ker\bar e}{\im\bar f}
 =\dim_\C\frac{\ker\bar f}{\im\bar e},
\qquad 0\le m\le\rank_R M.
\end{equation}
The module \(\Lambda\) is indecomposable and is simple exactly
when \(m=0\). The sequence \eqref{eq:non-weight-extension}
is nonsplit when \(m>0\).
\end{enumerate}
\end{theorem}
\begin{proof}
By Proposition~\ref{prop:evaluation},
\(\Der_\delta(\s,M)\cong\Lambda\ne0\).
For \(\delta\ne2\), Theorem~\ref{thm:translation}(i),(ii)
gives simplicity, since \(M\not\cong V(1)\).
The \(R\)-module \(\Lambda\subseteq M\) is torsion-free, so it has
no nonzero weight vectors. This proves (i).

For (ii), let \(t=0\), so \(\Omega\) acts by \(-1\) on \(M\). Equations
\eqref{eq:lattice} and \eqref{eq:transport} become
\begin{equation}\label{eq:two-lattice}
\begin{gathered}
\Lambda=\{u\in M:eu\in(h-1)M,\ fu\in(h+1)M\},\\
h*u=hu,\qquad e*u=\frac{h+1}{h-1}eu,\qquad
f*u=\frac{h-1}{h+1}fu.
\end{gathered}
\end{equation}
By the proof of Theorem~\ref{thm:translation}(iii), the operator
\(\mathrm J\) on \(\Hom_\C(\s,M)\) satisfies
\[
\im\mathrm J=\soc\Der_2(\s,M),
\]
and this socle is simple and essential. Its image under evaluation
at \(h\) can be computed from
\(\mathrm J=-2\mathrm S^2+3\mathrm S-1\) in \eqref{eq:nilpotent-operator}:
for \(F\in\Hom_\C(\s,M)\),
\[
\begin{aligned}
(\mathrm JF)(h)
&=(fh-3f)F(e)-(fe+ef+1)F(h)+(eh+3e)F(f)\\
&=(h-1)fF(e)+\frac{h^2-1}{2}F(h)+(h+1)eF(f),
\end{aligned}
\]
where the second equality uses \(2(ef+fe)+h^2=-1\) on \(M\).
The three values \(F(e),F(h),F(f)\) vary independently over \(M\).
Thus evaluation maps \(\im\mathrm J\) onto \(\Lambda_0\),
proving the socle assertion.

To compute \(\Lambda/\Lambda_0\), consider the maps in
\eqref{eq:fibers}. They are well defined because the shift
relations give
\[
e((h+1)M)\subseteq(h-1)M,\qquad
f((h-1)M)\subseteq(h+1)M.
\]
The identities
\[
fe=-\frac{1}{4}(h+1)^2,\qquad ef=-\frac{1}{4}(h-1)^2
\]
give \(\bar f\bar e=\bar e\bar f=0\).
Put \(P=(h^2-1)M\). Under the Chinese remainder isomorphism
\[
M/P\cong M/(h+1)M\oplus M/(h-1)M,
\]
the two conditions in \eqref{eq:two-lattice} restrict the first
coordinate to \(\ker\bar e\) and the second to \(\ker\bar f\).
For \(v\in M\), the generators of \(\Lambda_0/P\) have images
\[
\begin{aligned}
(h-1)fv+P&\longmapsto\bigl(-2(fv+(h+1)M),0\bigr),\\
(h+1)ev+P&\longmapsto\bigl(0,2(ev+(h-1)M)\bigr).
\end{aligned}
\]
Together with \eqref{eq:Lambda0}, this gives
\[
\Lambda/P\cong\ker\bar e\oplus\ker\bar f,\qquad
\Lambda_0/P\cong\im\bar f\oplus\im\bar e.
\]
Since \(h*u=hu\), the quotient \(Q=\Lambda/\Lambda_0\) has
weight spaces
\[
Q_{-1}\cong\frac{\ker\bar e}{\im\bar f},\qquad
Q_1\cong\frac{\ker\bar f}{\im\bar e},\qquad
Q_\nu=0\quad(\nu\ne\pm1).
\]
By Theorem~\ref{thm:translation}(iii) and the socle assertion above,
\(Q\cong V(1)^{\oplus m}\) for some finite \(m\).
Thus \(m=\dim Q_{-1}=\dim Q_1\), proving
\eqref{eq:non-weight-extension} and both expressions for \(m\).

To bound \(m\), fix \(\alpha\in\C\) and choose
\(v_1,\ldots,v_r\in M\) whose classes in
\(M/(h-\alpha)M\) are linearly independent.
If these lifts were \(K\)-linearly dependent, clearing denominators
would give
\[
\sum_{i=1}^r a_i(h)v_i=0,\qquad a_i(h)\in R,
\]
with the coefficients not all zero. Since \(M\) is torsion-free,
we may cancel their largest common power of \(h-\alpha\) and
assume that at least one \(a_i(\alpha)\) is nonzero. Reduction modulo
\((h-\alpha)M\) then gives
\[
\sum_{i=1}^r a_i(\alpha)\bigl(v_i+(h-\alpha)M\bigr)=0,
\]
contrary to the choice of the classes. The lifts are therefore
\(K\)-linearly independent, so
\(\dim_\C M/(h-\alpha)M\le\rank_R M\). In particular,
\[
m\le\dim_\C M/(h+1)M\le\rank_R M,
\]
which proves the bound in \eqref{eq:multiplicity}.
Every nonzero submodule of \(\Lambda\) meets the simple socle
\(\Lambda_0\), so \(\Lambda\) is indecomposable. If \(m=0\),
then \(\Lambda=\Lambda_0\) is simple; if \(m>0\), then
\(\Lambda_0\) is proper, and its essentiality rules out a
splitting of \eqref{eq:non-weight-extension}.
\end{proof}

\subsection{Block's classification}\label{subsec:block}

To identify \(\Lambda\) and its socle in Block's classification
\cite{Bavula90,Block81}, we use the notation of \cite{Maz10}.
Set
\begin{equation}\label{eq:B}
B=K[X,X^{-1};\sigma],\qquad \sigma(h)=h-2,\qquad
q_d(h)=\frac{d-h^2-2h}{4}\quad(d\in\C).
\end{equation}
Thus \(Xr(h)=r(h-2)X\). The embedding
\begin{equation}\label{eq:embedding}
A_d\hookrightarrow B,\qquad
e\longmapsto X,\qquad f\longmapsto q_d(h)X^{-1},\qquad h\longmapsto h
\end{equation}
identifies \(B\) with the localization of \(A_d\) at
\(R\setminus\{0\}\), since \(X^{-1}=q_d(h)^{-1}f\).
With our Casimir normalization, the parameter in \cite{Maz10}
satisfies \((\lambda+1)^2=d+1\).
For an irreducible nonunit \(b\in B\), put
\begin{equation}\label{eq:BB}
W_b=B/Bb,\qquad L_d(b)=\soc_{A_d}W_b.
\end{equation}
For each \(d\in\C\), we regard \(W_b\) as an \(A_d\)-module
through the embedding \eqref{eq:embedding}.
By \cite[Proposition~6.2.8 and Theorem~6.3.6]{Maz10},
\(L_d(b)\) is simple, \(K\otimes_R L_d(b)=W_b\), and every
simple non-weight \(A_d\)-module is obtained in this way.
The socle is essential. Indeed, for each nonzero \(v\in W_b\),
the equality \(K\otimes_R L_d(b)=W_b\) gives a nonzero
\(r(h)\in R\) such that \(r(h)v\in L_d(b)\).
Since \(W_b\) is a \(K\)-vector space, \(r(h)v\ne0\).
Thus every nonzero \(A_d\)-submodule of \(W_b\) meets
\(L_d(b)\).
Write \(b=\sum_{r=i}^j a_r(h)X^r\), with \(i<j\) and
\(a_i a_j\ne0\). For every \(\ell\in\Z\), we have
\[
0=X^\ell b(1+Bb)
=\sum_{r=i}^j a_r(h-2\ell)\bigl(X^{r+\ell}+Bb\bigr).
\]
Since the two end coefficients are invertible in \(K\), recursion
in both directions expresses every \(X^k+Bb\) in the \(K\)-span of
\(X^i+Bb,\ldots,X^{j-1}+Bb\). Hence
\[
\rank_R L_d(b)=\dim_K W_b\le j-i<\infty.
\]

\begin{corollary}\label{cor:block-modules}
Let \(t=2/\delta-1\), let \(M=L_{t^2-1}(b)\subset W_b\),
and let \(\Lambda\) be given by \eqref{eq:lattice}.
The \(\delta\)-derivations with values in \(M\) are parametrized
bijectively by
\begin{equation}\label{eq:BB-parameters}
w\in(h-t+1)^{-1}\Lambda\subset W_b
\end{equation}
via the formulas
\begin{equation}\label{eq:BB-D}
\begin{aligned}
D_w(e)&=Xw,\\
D_w(h)&=(h-t+1)w,\\
D_w(f)&=-\frac{1}{4}(h-t+3)(h-t+1)X^{-1}w.
\end{aligned}
\end{equation}
Equivalently, the three values in \eqref{eq:BB-D} must lie in \(M\).
\begin{enumerate}
\item If \(\delta\ne2\), the parameter \(w\) ranges over
\(L_{(t-2)^2-1}(b)\), and
\begin{equation}\label{eq:BB-simple}
\Der_\delta(\s,L_{t^2-1}(b))\cong L_{(t-2)^2-1}(b).
\end{equation}
\item If \(\delta=2\), there is an exact sequence
\begin{equation}\label{eq:BB-extension}
0\longrightarrow L_3(b)\longrightarrow
\Der_2(\s,L_{-1}(b))\longrightarrow V(1)^{\oplus m_b}
\longrightarrow0,
\end{equation}
where \(m_b\) is the multiplicity in \eqref{eq:multiplicity}
for \(M=L_{-1}(b)\). The inclusion is \(w\mapsto D_w\),
and its image is the simple essential socle.
\end{enumerate}
\end{corollary}
\begin{proof}
Write \(c=t^2-1\) and \(c'=(t-2)^2-1\), as in \eqref{eq:parameters}.
Substituting \(u=(h-t+1)w\) in \eqref{eq:Du} and using
\[
e((h-t+1)w)=(h-t-1)Xw,\qquad
q_c(h)=-\frac{1}{4}(h-t+1)(h+t+1)
\]
gives \eqref{eq:BB-D}. For the \(f\)-action in
\eqref{eq:transport}, we obtain
\[
\begin{aligned}
&(h-t+1)^{-1}\frac{h+t-1}{h+t+1}q_c(h)X^{-1}(h-t+1)\\
&\qquad=\frac{(h+t-1)(h-t+3)}{(h-t+1)(h+t+1)}q_c(h)X^{-1}\\
&\qquad=q_{c'}(h)X^{-1}.
\end{aligned}
\]
Together with the \(h\)- and \(e\)-actions, this gives
\begin{equation}\label{eq:BB-action}
h*w=hw,\qquad e*w=Xw,\qquad f*w=q_{c'}(h)X^{-1}w.
\end{equation}
This is the \(A_{c'}\)-module structure on \(W_b\).
If \(\delta\ne2\), Theorem~\ref{thm:non-weight}(i) shows that
\((h-t+1)^{-1}\Lambda\) is a simple \(A_{c'}\)-submodule of \(W_b\).
Since \(L_{c'}(b)\) is the simple essential socle of \(W_b\),
it follows that \((h-t+1)^{-1}\Lambda=L_{c'}(b)\).
If \(\delta=2\), the change of variables \(u=(h+1)w\)
makes \((h+1)^{-1}\Lambda_0\) a simple \(A_3\)-submodule
of \(W_b\). It therefore equals \(L_3(b)\), so
\begin{equation}\label{eq:socle-coordinates}
\Lambda_0=(h+1)L_3(b).
\end{equation}
Now \eqref{eq:non-weight-extension} gives \eqref{eq:BB-extension}.
\end{proof}

\begin{remark}
For \(M=L_{-1}(b)\), the exact sequences
\eqref{eq:translated-exact} and \eqref{eq:BB-extension} give
\[
[\T(M):L_3(b)]=2,\qquad [\T(M):V(1)]=m_b.
\]
\end{remark}

For \(M=L_{-1}(b)\) of \(R\)-rank one, \eqref{eq:multiplicity} gives
\(m_b\in\{0,1\}\). The following examples realize both values.

\begin{example}\label{ex:whittaker}
For \(\eta\in\C^\times\) and \(\kappa\in\C\), the
simple Whittaker module \(W_{\eta,\kappa}=\C[h]\)
\cite{Maz10} has action
\begin{equation}\label{eq:whittaker-source}
h\cdot p=hp,\qquad e\cdot p=\eta p(h-2),\qquad
f\cdot p=\frac{\kappa-h^2-2h}{4\eta}p(h+2).
\end{equation}
This module is non-weight and \(R\)-free of rank one, with
\(e\cdot1=\eta\) and \(\Omega=\kappa\).
The parameters in the cited model are \(\alpha=\eta^{-1}\) and
\((\lambda+1)^2=\kappa+1\).

By Proposition~\ref{prop:evaluation}, nonzero derivations exist
exactly when \(\kappa=t^2-1\). In this case, the condition
\(fu\in(h+t+1)W_{\eta,\kappa}\) is automatic, whereas
\(eu\in(h-t-1)W_{\eta,\kappa}\) is equivalent to
\(u\in(h-t+1)\C[h]\). Writing \(u=(h-t+1)p(h)\), we obtain
all derivations in the form
\begin{equation}\label{eq:whittaker-derivations}
\begin{aligned}
D_p(e)&=\eta p(h-2),\\
D_p(h)&=(h-t+1)p(h),\\
D_p(f)&=-\frac{(h-t+1)(h-t+3)}{4\eta}p(h+2),
\end{aligned}
\qquad p\in\C[h].
\end{equation}
Under this parametrization, \eqref{eq:transport} gives
\[
h*p=hp,\qquad e*p=\eta p(h-2),\qquad
f*p=\frac{(t-2)^2-1-h^2-2h}{4\eta}p(h+2).
\]
Consequently
\begin{equation}\label{eq:whittaker-target}
\Der_\delta(\s,W_{\eta,t^2-1})
\cong W_{\eta,(t-2)^2-1}.
\end{equation}
At \(\delta=2\), the target is the simple Whittaker module
\(W_{\eta,3}\). Hence the quotient in
\eqref{eq:non-weight-extension} is zero.
\end{example}

\begin{example}\label{ex:nonsimple}
On \(M=\C[h]\), put
\begin{equation}\label{eq:rank-one-source}
h\cdot p=hp,\qquad e\cdot p=\eta(h-1)p(h-2),\qquad
f\cdot p=-\frac{h+1}{4\eta}p(h+2),\qquad \eta\ne0.
\end{equation}
These operators satisfy the relations in \eqref{eq:sl2}, and
\(\Omega\) acts by \(-1\).
Since \(h\) acts by multiplication, every nonzero submodule is
an ideal \((p(h))\) of \(\C[h]\), with \(p\ne0\).
Suppose \(p\) is nonconstant and let \(\alpha\) be a root.
Stability under the action forces \(p(\alpha-2)=0\) unless \(\alpha=1\),
and \(p(\alpha+2)=0\) unless \(\alpha=-1\).
The descending sequence of roots must terminate at \(1\),
and the ascending sequence at \(-1\), since \(p\) has only
finitely many roots. Thus
\(\alpha=1+2k=-1-2\ell\) for some \(k,\ell\in\Z_{\ge0}\),
a contradiction.
Hence \(p\) is a nonzero constant and \(M\) is simple.

Here \(\Lambda=\C[h]\), so all \(2\)-derivations are
\[
D_p(e)=\eta p(h-2),\qquad D_p(h)=p(h),\qquad
D_p(f)=-\frac{1}{4\eta}p(h+2).
\]
The induced action on the parameter \(p\) is
\begin{equation}\label{eq:rank-one-target}
h*p=hp,\qquad e*p=\eta(h+1)p(h-2),\qquad
f*p=-\frac{h-1}{4\eta}p(h+2).
\end{equation}
Formula \eqref{eq:Lambda0} gives
\(\Lambda_0=(h^2-1)\C[h]\). Both spaces in \eqref{eq:fibers} have dimension one, and
\(\bar e=\bar f=0\). Hence \eqref{eq:multiplicity} gives \(m=1\).
Let \(v_+,v_-\) be the standard basis of \(V(1)\), with
\(hv_\pm=\pm v_\pm\), \(ev_-=v_+\), and \(fv_+=v_-\).
The map
\[
\Der_2(\s,M)\longrightarrow V(1),\qquad
D_p\longmapsto p(1)v_++2\eta p(-1)v_-
\]
is surjective. By \eqref{eq:rank-one-target},
\[
\begin{aligned}
(e*p)(1)&=2\eta p(-1),& (e*p)(-1)&=0,\\
(f*p)(1)&=0,& (f*p)(-1)&=\frac{p(1)}{2\eta}.
\end{aligned}
\]
Together with \(h*p=hp\), these identities show that the map is
\(\s\)-equivariant. Its kernel consists of the \(D_p\) with
\(p(1)=p(-1)=0\), or equivalently \(p\in(h^2-1)\C[h]\).
Thus we obtain the nonsplit sequence of Theorem~\ref{thm:non-weight},
with first map \(p\mapsto D_p\):
\[
0\longrightarrow(h^2-1)\C[h]\longrightarrow
\Der_2(\s,M)\longrightarrow V(1)\longrightarrow0.
\]
\end{example}

\end{document}